# SPLINE-BACKFITTED KERNEL SMOOTHING OF NONLINEAR ADDITIVE AUTOREGRESSION MODEL[1]

BY LI WANG AND LIJIAN YANG

*University of Georgia and Michigan State University*

Application of nonparametric and semiparametric regression techniques to high-dimensional time series data has been hampered due to the lack of effective tools to address the "curse of dimensionality." Under rather weak conditions, we propose spline-backfitted kernel estimators of the component functions for the nonlinear additive time series data that are both computationally expedient so they are usable for analyzing very high-dimensional time series, and theoretically reliable so inference can be made on the component functions with confidence. Simulation experiments have provided strong evidence that corroborates the asymptotic theory.

**1. Introduction.** For the past three decades, various nonparametric and semiparametric regression techniques have been developed for the analysis of nonlinear time series; see, for example, [14, 21, 25], to name one article representative of each decade. Application to high-dimensional time series data, however, has been hampered due to the scarcity of smoothing tools that are not only computationally expedient but also theoretically reliable, which has motivated the proposed procedures of this paper.

In high-dimensional time series smoothing, one unavoidable issue is the "curse of dimensionality," which refers to the poor convergence rate of nonparametric estimation of general multivariate functions. One solution is regression in the form of an additive model introduced by [9]:

$$Y_i = m(X_{i1}, \ldots, X_{id}) + \sigma(X_{i1}, \ldots, X_{id})\varepsilon_i,$$

(1.1)

$$m(x_1, \ldots, x_d) = c + \sum_{\alpha=1}^{d} m_\alpha(x_\alpha),$$

Received December 2005; revised January 2007.
[1]Supported in part by NSF Grant DMS-04-05330.

in which the sequence $\{Y_i, \mathbf{X}_i^T\}_{i=1}^n = \{Y_i, X_{i1}, \ldots, X_{id}\}_{i=1}^n$ is a length-$n$ realization of a $(d+1)$-dimensional time series, the $d$-variate functions $m$ and $\sigma$ are the mean and standard deviation of the response $Y_i$ conditional on the predictor vector $\mathbf{X}_i = \{X_{i1}, \ldots, X_{id}\}^T$, and each $\varepsilon_i$ is a white noise conditional on $\mathbf{X}_i$. In a nonlinear additive autoregression data-analytical context, each predictor $X_{i\alpha}, 1 \leq \alpha \leq d$, could be observed lagged values of $Y_i$, such as $X_{i\alpha} = Y_{i-\alpha}$, or of an exogenous time series. Model (1.1), therefore, is the exact same nonlinear additive autoregression model of [14] and [2] with exogenous variables. For identifiability, additive component functions must satisfy the conditions $Em_\alpha(X_{i\alpha}) \equiv 0$, $\alpha = 1, \ldots, d$.

We propose estimators of the unknown component functions $\{m_\alpha(\cdot)\}_{\alpha=1}^d$ based on a geometrically $\alpha$-mixing sample $\{Y_i, X_{i1}, \ldots, X_{id}\}_{i=1}^n$ following model (1.1). If the data were actually i.i.d. observations instead of a time series realization, many methods would be available for estimating $\{m_\alpha(\cdot)\}_{\alpha=1}^d$. For instance, there are four types of kernel-based estimators: the classic backfitting estimators (CBE) of [9] and [19]; marginal integration estimators (MIE) of [6, 16, 17, 22, 30] and a kernel-based method of estimating rate to optimality of [10]; the smoothing backfitting estimators (SBE) of [18]; and the two-stage estimators, such as one step backfitting of the integration estimators of [15], one step backfitting of the projection estimators of [11] and one Newton step from the nonlinear LSE estimators of [12]. For the spline estimators, see [13, 23, 24] and [28].

In the time series context, however, there are fewer theoretically justified methods due to the additional difficulty posed by dependence in the data. Some of these are the kernel estimators via marginal integration of [25, 29], and the spline estimators of [14]. In addition, [27] has extended the marginal integration kernel estimator to additive coefficient models for weakly dependent data. All of these existing methods are unsatisfactory in regard to either the computational or the theoretical issue. The existing kernel methods are too computationally intensive for high dimension $d$, thus limiting their applicability to a small number of predictors. Spline methods, on the other hand, provide only convergence rates but no asymptotic distributions, so no measures of confidence can be assigned to the estimators.

If the last $d-1$ component functions were known by "oracle," one could create $\{Y_{i1}, X_{i1}\}_{i=1}^n$ with $Y_{i1} = Y_i - c - \sum_{\alpha=2}^d m_\alpha(X_{i\alpha}) = m_1(X_{i1}) + \sigma(X_{i1}, \ldots, X_{id})\varepsilon_i$, from which one could compute an "oracle smoother" to estimate the only unknown function $m_1(x_1)$, thus effectively bypassing the "curse of dimensionality." The idea of [15] was to obtain an approximation to the unobservable variables $Y_{i1}$ by substituting $m_\alpha(X_{i\alpha})$, $i = 1, \ldots, n$, $\alpha = 2, \ldots, d$, with marginal integration kernel estimates and arguing that the error incurred by this "cheating" is of smaller magnitude than the rate $O(n^{-2/5})$ for estimating the function $m_1(x_1)$ from the unobservable data. We modify the procedure of [15] by substituting $m_\alpha(X_{i\alpha})$, $i = 1, \ldots, n$, $\alpha = 2, \ldots, d$, with spline estimators. Specifically, we propose a two-stage estimation procedure; first we pre-estimate $\{m_\alpha(x_\alpha)\}_{\alpha=2}^d$ by its pilot estimator through an undersmoothed centered standard spline procedure; next we construct the pseudo response $\hat{Y}_{i1}$ and approximate $m_1(x_1)$ by its Nadaraya–Watson estimator in (2.12).

The above proposed spline-backfitted kernel (SPBK) estimation method has several advantages compared to most of the existing methods. First, as pointed out in [22], the estimator of [15] mixed up different projections, making it uninterpretable if the real data generating process deviates from additivity, while the projections in both steps of our estimator are with respect to the same measure. Second, since our pilot spline estimator is thousands of times faster than the pilot kernel estimator in [15], our proposed method is computationally expedient; see Table 2. Third, the SPBK estimator can be shown to be as efficient as the "oracle smoother" uniformly over any compact range, whereas [15] proved such "oracle efficiency" only at a single point. Moreover, the regularity conditions in our paper are natural and appealing and close to being minimal. In contrast, higher-order smoothness is needed with growing dimensionality of the regressors in [17]. Stronger and more obscure conditions are assumed for the two-stage estimation proposed by [12].

The SPBK estimator achieves its seemingly surprising success by borrowing the strengths of both spline and kernel: the spline does a quick initial estimation of all additive components and removes them all except the one of interest; kernel smoothing is then applied to the cleaned univariate data to estimate with asymptotic distribution. Propositions 4.1 and 5.1 are the keys in understanding the proposed estimators' uniform oracle efficiency. They accomplish the well-known "reducing bias by undersmoothing" in the first step using spline and "averaging out the variance" in the second step with kernel, both steps taking advantage of the joint asymptotics of kernel and spline functions, which is the new feature of our proofs.

Reference [7] provides generalized likelihood ratio (GLR) tests for additive models using the backfitting estimator. A similar GLR test based on our SPBK estimator is feasible for future research.

The rest of the paper is organized as follows. In Section 2 we introduce the SPBK estimator and state its asymptotic "oracle efficiency" under appropriate assumptions. In Section 3 we provide some insights into the ideas behind our proofs of the main results, by decomposing the estimator's "cheating" error into a bias and a variance part. In Section 4 we show the uniform order of the bias term. In Section 5 we show the uniform order of the variance term. In Section 6 we present Monte Carlo results to demonstrate that the SPBK estimator does indeed possess the claimed asymptotic properties. All technical proofs are contained in the Appendix.

**2. The SPBK estimator.** In this section we describe the spline-backfitted kernel estimation procedure. For convenience, we denote vectors as $\mathbf{x} = (x_1,\ldots,x_d)$ and take $\|\cdot\|$ as the usual Euclidean norm on $R^d$, that is, $\|\mathbf{x}\| = \sqrt{\sum_{\alpha=1}^d x_\alpha^2}$, and $\|\cdot\|_\infty$ the sup norm, that is, $\|\mathbf{x}\|_\infty = \sup_{1\le\alpha\le d}|x_\alpha|$. In what follows, let $Y_i$ and $\mathbf{X}_i = (X_{i1},\ldots,X_{id})^T$ be the $i$th response and predictor vector. Denote by $\mathbf{Y} = (Y_1,\ldots,Y_n)^T$ the response vector and $(\mathbf{X}_1,\ldots,\mathbf{X}_n)^T$ the design matrix.

Let $\{Y_i,\mathbf{X}_i^T\}_{i=1}^n = \{Y_i,X_{i1},\ldots,X_{id}\}_{i=1}^n$ be observations from a geometrically $\alpha$-mixing process following model (1.1). We assume that the predictor $X_\alpha$ is distributed on a compact interval $[a_\alpha,b_\alpha]$, $\alpha = 1,\ldots,d$, and without loss of generality, we take all intervals $[a_\alpha,b_\alpha] = [0,1]$, $\alpha = 1,\ldots,d$. We preselect an integer $N = N_n \sim n^{2/5}\log n$; see assumption (A6) below. Next, we define for any $\alpha = 1,\ldots,d$ the first-order B spline function ([3], page 89), or say the constant B spline function is the indicator function $I_{J,\alpha}(x_\alpha)$ of the $N+1$ equally spaced subintervals of the finite interval $[0,1]$ with length $H = H_n = (N+1)^{-1}$, that is,

$$(2.1) \quad I_{J,\alpha}(x_\alpha) = \begin{cases} 1, & JH \le x_\alpha < (J+1)H, \\ 0, & \text{otherwise,} \end{cases} \quad J = 0,1,\ldots,N.$$

Define the following centered spline basis:

$$(2.2) \quad b_{J,\alpha}(x_\alpha) = I_{J+1,\alpha}(x_\alpha) - \frac{\|I_{J+1,\alpha}\|_2}{\|I_{J,\alpha}\|_2} I_{J,\alpha}(x_\alpha)$$

$$\forall \alpha = 1,\ldots,d,\ J = 1,\ldots,N,$$

with the standardized version given for any $\alpha = 1,\ldots,d$,

$$(2.3) \quad B_{J,\alpha}(x_\alpha) = \frac{b_{J,\alpha}(x_\alpha)}{\|b_{J,\alpha}\|_2} \quad \forall J = 1,\ldots,N.$$

Define next the $(1+dN)$-dimensional space $G = G[0,1]$ of additive spline functions as the linear space spanned by $\{1,B_{J,\alpha}(x_\alpha),\alpha=1,\ldots,d,J=1,\ldots,N\}$, and denote by $G_n \subset R^n$ the linear space spanned by $\{1,\{B_{J,\alpha}(X_{i\alpha})\}_{i=1}^n, \alpha=1,\ldots,d, J=1,\ldots,N\}$. As $n \to \infty$, the dimension of $G_n$ becomes $1+dN$ with probability approaching 1. The spline estimator of the additive function $m(\mathbf{x})$ is the unique element $\hat{m}(\mathbf{x}) = \hat{m}_n(\mathbf{x})$ from the space $G$ so that the vector $\{\hat{m}(\mathbf{X}_1),\ldots,\hat{m}(\mathbf{X}_n)\}^T$ best approximates the response vector $\mathbf{Y}$. To be precise, we define

$$(2.4) \quad \hat{m}(\mathbf{x}) = \hat{\lambda}'_0 + \sum_{\alpha=1}^d \sum_{J=1}^N \hat{\lambda}'_{J,\alpha} I_{J,\alpha}(x_\alpha),$$

where the coefficients $(\hat{\lambda}'_0, \hat{\lambda}'_{1,1}, \ldots, \hat{\lambda}'_{N,d})$ are solutions of the least squares problem

$$\{\hat{\lambda}'_0, \hat{\lambda}'_{1,1}, \ldots, \hat{\lambda}'_{N,d}\}^T = \underset{R^{dN+1}}{\arg\min} \sum_{i=1}^{n} \left\{Y_i - \lambda_0 - \sum_{\alpha=1}^{d} \sum_{J=1}^{N} \lambda_{J,\alpha} I_{J,\alpha}(X_{i\alpha})\right\}^2.$$

Simple linear algebra shows that

$$(2.5) \qquad \hat{m}(\mathbf{x}) = \hat{\lambda}_0 + \sum_{\alpha=1}^{d} \sum_{J=1}^{N} \hat{\lambda}_{J,\alpha} B_{J,\alpha}(x_\alpha),$$

where $(\hat{\lambda}_0, \hat{\lambda}_{1,1}, \ldots, \hat{\lambda}_{N,d})$ are solutions of the least squares problem

$$\{\hat{\lambda}_0, \hat{\lambda}_{1,1}, \ldots, \hat{\lambda}_{N,d}\}^T = \underset{R^{dN+1}}{\arg\min} \sum_{i=1}^{n} \left\{Y_i - \lambda_0 - \sum_{\alpha=1}^{d} \sum_{J=1}^{N} \lambda_{J,\alpha} B_{J,\alpha}(X_{i\alpha})\right\}^2;$$

(2.6)

while (2.4) is used for data-analytic implementation, the mathematically equivalent expression (2.5) is convenient for asymptotic analysis.

The pilot estimators of each component function and the constant are

$$(2.7) \qquad \hat{m}_\alpha(x_\alpha) = \sum_{J=1}^{N} \hat{\lambda}_{J,\alpha} B_{J,\alpha}(x_\alpha) - n^{-1} \sum_{i=1}^{n} \sum_{J=1}^{N} \hat{\lambda}_{J,\alpha} B_{J,\alpha}(X_{i\alpha}),$$

$$\hat{m}_c = \hat{\lambda}_0 + n^{-1} \sum_{\alpha=1}^{d} \sum_{i=1}^{n} \sum_{J=1}^{N} \hat{\lambda}_{J,\alpha} B_{J,\alpha}(X_{i\alpha}).$$

These pilot estimators are then used to define new pseudo-responses $\hat{Y}_{i1}$, which are estimates of the unobservable "oracle" responses $Y_{i1}$. Specifically,

$$(2.8) \quad \hat{Y}_{i1} = Y_i - \hat{c} - \sum_{\alpha=2}^{d} \hat{m}_\alpha(X_{i\alpha}), \qquad Y_{i1} = Y_i - c - \sum_{\alpha=2}^{d} m_\alpha(X_{i\alpha}),$$

where $\hat{c} = \overline{Y}_n = n^{-1} \sum_{i=1}^{n} Y_i$, which is a $\sqrt{n}$-consistent estimator of $c$ by the central limit theorem. Next, we define the spline-backfitted kernel estimator of $m_1(x_1)$ as $\hat{m}_1^*(x_1)$ based on $\{\hat{Y}_{i1}, X_{i1}\}_{i=1}^{n}$, which attempts to mimic the would-be Nadaraya–Watson estimator $\tilde{m}_1^*(x_1)$ of $m_1(x_1)$ based on $\{Y_{i1}, X_{i1}\}_{i=1}^{n}$ if the unobservable "oracle" responses $\{Y_{i1}\}_{i=1}^{n}$ were available:

$$\hat{m}_1^*(x_1) = \frac{\sum_{i=1}^{n} K_h(X_{i1} - x_1)\hat{Y}_{i1}}{\sum_{i=1}^{n} K_h(X_{i1} - x_1)},$$

(2.9)

$$\tilde{m}_1^*(x_1) = \frac{\sum_{i=1}^{n} K_h(X_{i1} - x_1)Y_{i1}}{\sum_{i=1}^{n} K_h(X_{i1} - x_1)},$$

where $\hat{Y}_{i1}$ and $Y_{i1}$ are defined in (2.8).

Throughout this paper, on any fixed interval $[a,b]$, we denote the space of second-order smooth functions as $C^{(2)}[a,b] = \{m|m'' \in C[a,b]\}$ and the class of Lipschitz continuous functions for any fixed constant $C > 0$ as $\text{Lip}([a,b], C) = \{m||m(x) - m(x')| \leq C|x - x'|, \; \forall x, x' \in [a,b]\}$.

Before presenting the main results, we state the following assumptions.

(A1) The additive component function $m_1(x_1) \in C^{(2)}[0,1]$, while there is a constant $0 < C_\infty < \infty$ such that $m_\beta \in \text{Lip}([0,1], C_\infty), \; \forall \beta = 2, \ldots, d$.

(A2) There exist positive constants $K_0$ and $\lambda_0$ such that $\alpha(n) \leq K_0 e^{-\lambda_0 n}$ holds for all $n$, with the $\alpha$-mixing coefficients for $\{\mathbf{Z}_i = (\mathbf{X}_i^T, \varepsilon_i)\}_{i=1}^n$ defined as

$$(2.10) \quad \alpha(k) = \sup_{B \in \sigma\{\mathbf{Z}_s, s \leq t\}, C \in \sigma\{\mathbf{Z}_s, s \geq t+k\}} |P(B \cap C) - P(B)P(C)|, \quad k \geq 1.$$

(A3) The noise $\varepsilon_i$ satisfies $E(\varepsilon_i|\mathbf{X}_i) = 0, E(\varepsilon_i^2|\mathbf{X}_i) = 1$ and $E(|\varepsilon_i|^{2+\delta}|\mathbf{X}_i) < M_\delta$ for some $\delta > 1/2$ and a finite positive $M_\delta$ and $\sigma(\mathbf{x})$ is continuous on $[0,1]^d$:

$$0 < c_\sigma \leq \inf_{\mathbf{x} \in [0,1]^d} \sigma(\mathbf{x}) \leq \sup_{\mathbf{x} \in [0,1]^d} \sigma(\mathbf{x}) \leq C_\sigma < \infty.$$

(A4) The density function $f(\mathbf{x})$ of $\mathbf{X}$ is continuous and

$$0 < c_f \leq \inf_{\mathbf{x} \in [0,1]^d} f(\mathbf{x}) \leq \sup_{\mathbf{x} \in [0,1]^d} f(\mathbf{x}) \leq C_f < \infty.$$

The marginal densities $f_\alpha(x_\alpha)$ of $X_\alpha$ have continuous derivatives on $[0,1]$ as well as the uniform upper bound $C_f$ and lower bound $c_f$.

(A5) The kernel function $K \in \text{Lip}([-1,1], C_\infty)$ for some constant $C_k > 0$, and is bounded, nonnegative, symmetric and supported on $[-1,1]$. The bandwidth $h \sim n^{-1/5}$, that is, $c_h n^{-1/5} \leq h \leq C_h n^{-1/5}$ for some positive constants $C_h, c_h$.

(A6) The number of interior knots $N \sim n^{2/5} \log n$, that is, $c_N n^{2/5} \log n \leq N \leq C_N n^{2/5} \log n$ for some positive constants $c_N, C_N$.

REMARK 2.1. The smoothness assumption of the true component functions is greatly relaxed in our paper and we believe that our assumption (A1) is close to being minimal. By the result of [20], a geometrically ergodic time series is a strongly mixing sequence. Therefore, assumption (A2) is suitable for (1.1) as a time series model under the aforementioned assumptions. Assumptions (A3)–(A5) are typical in the nonparametric smoothing literature; see, for instance, [5]. For (A6), the proof of Theorem 2.1 in the Appendix will make it clear that the number of knots can be of the more general form $N \sim n^{2/5} N'$, where the sequence $N'$ satisfies $N' \to \infty$, $n^{-\theta} N' \to 0$ for any $\theta > 0$. There is no optimal way to choose $N'$ as in the literature. Here we select $N$ to be of barely larger order than $n^{2/5}$.

The asymptotic property of the kernel smoother $\tilde{m}_1^*(x_1)$ is well developed. Under assumptions (A1)–(A5), it is straightforward to verify (as in [1]) that

$$\sup_{x_1 \in [h, 1-h]} |\tilde{m}_1^*(x_1) - m_1(x_1)| = o_p(n^{-2/5} \log n),$$

$$\sqrt{nh}\{\tilde{m}_1^*(x_1) - m_1(x_1) - b_1(x_1)h^2\} \xrightarrow{D} N\{0, v_1^2(x_1)\},$$

where

(2.11)
$$b_1(x_1) = \int u^2 K(u)\, du \{m_1''(x_1)f_1(x_1)/2 + m_1'(x_1)f_1'(x_1)\} f_1^{-1}(x_1),$$
$$v_1^2(x_1) = \int K^2(u)\, du\, E[\sigma^2(X_1, \ldots, X_d)|X_1 = x_1] f_1^{-1}(x_1).$$

The following theorem states that the asymptotic uniform magnitude of the difference between $\hat{m}_1^*(x_1)$ and $\tilde{m}_1^*(x_1)$ is of order $o_p(n^{-2/5})$, which is dominated by the asymptotic uniform size of $\tilde{m}_1^*(x_1) - m_1(x_1)$. As a result, $\hat{m}_1^*(x_1)$ will have the same asymptotic distribution as $\tilde{m}_1^*(x_1)$.

THEOREM 2.1. *Under assumptions* (A1)–(A6), *the SPBK estimator* $\hat{m}_1^*(x_1)$ *given in* (2.9) *satisfies*

$$\sup_{x_1 \in [0,1]} |\hat{m}_1^*(x_1) - \tilde{m}_1^*(x_1)| = o_p(n^{-2/5}).$$

*Hence with* $b_1(x_1)$ *and* $v_1^2(x_1)$ *as defined in* (2.11), *for any* $x_1 \in [h, 1-h]$

$$\sqrt{nh}\{\hat{m}_1^*(x_1) - m_1(x_1) - b_1(x_1)h^2\} \xrightarrow{D} N\{0, v_1^2(x_1)\}.$$

REMARK 2.2. Theorem 2.1 holds for $\hat{m}_\alpha^*(x_\alpha)$ similarly constructed as $\hat{m}_1^*(x_1)$, for any $\alpha = 2, \ldots, d$, that is,

(2.12)
$$\hat{m}_\alpha^*(x_\alpha) = \frac{\sum_{i=1}^n K_h(X_{i\alpha} - x_\alpha)\hat{Y}_{i\alpha}}{\sum_{i=1}^n K_h(X_{i1} - x_\alpha)},$$
$$\hat{Y}_{i\alpha} = Y_i - \hat{c} - \sum_{1 \leq \beta \leq d, \beta \neq \alpha} \hat{m}_\beta(X_{i\beta}),$$

where $\hat{m}_\beta(X_{i\beta})$, $\beta = 1, \ldots, d$, are the pilot estimators of each component function given in (2.7). Similar constructions can be based on a local polynomial instead of the Nadaraya–Watson estimator. For more on the properties of local polynomial estimators, in particular, their minimax efficiency, see [5].

REMARK 2.3. Compared to the SBE in [18], the variance term $v_1(x_1)$ is identical to that of SBE and the bias term $b_1(x_1)$ is much more explicit than that of SBE, at least when the Nadaraya–Watson smoother is used. Theorem 2.1 can be used to construct asymptotic confidence intervals. Under assumptions (A1)–(A6), for any $\alpha \in (0,1)$, an asymptotic $100(1-\alpha)\%$ pointwise confidence interval for $m_1(x_1)$ is

$$(2.13) \quad \hat{m}_1^*(x_1) - b_1(x_1)h^2 \pm z_{\alpha/2}\hat{\sigma}_1(x_1)\left\{\int K^2(u)\,du\right\}^{1/2}\bigg/\{nh\hat{f}_1(x_1)\}^{1/2},$$

where $\hat{\sigma}_1(x_1)$ and $\hat{f}_1(x_1)$ are estimators of $E[\sigma^2(X_1,\ldots,X_d)|X_1 = x_1]$ and $f_1(x_1)$.

The following corollary provides the asymptotic distribution of $\hat{m}^*(\mathbf{x})$. The proof of this corollary is straightforward and therefore omitted.

COROLLARY 2.1. *Under assumptions* (A1)–(A6) *and the additional assumption that* $m_\alpha(x_\alpha) \in C^{(2)}[0,1], \alpha = 2,\ldots,d$, *for any* $\mathbf{x} \in [0,1]^d$, *the SPBK estimator* $\hat{m}_\alpha^*(\mathbf{x}), \alpha = 1,\ldots,d$, *is defined as given in* (2.12). *Let*

$$\hat{m}^*(\mathbf{x}) = \hat{c} + \sum_{\alpha=1}^d \hat{m}_\alpha^*(x_\alpha), \qquad b(\mathbf{x}) = \sum_{\alpha=1}^d b_\alpha(x_\alpha), \qquad v^2(\mathbf{x}) = \sum_{\alpha=1}^d v_\alpha^2(x_\alpha).$$

*Then*

$$\sqrt{nh}\{\hat{m}^*(\mathbf{x}) - m(\mathbf{x}) - b(\mathbf{x})h^2\} \xrightarrow{D} N\{0, v^2(\mathbf{x})\}.$$

**3. Decomposition.** In this section we introduce some additional notation to shed some light on the ideas behind the proof of Theorem 2.1. For any functions $\phi, \varphi$ on $[0,1]^d$, define the empirical inner product and the empirical norm as

$$\langle \phi, \varphi \rangle_{2,n} = n^{-1}\sum_{i=1}^n \phi(\mathbf{X}_i)\varphi(\mathbf{X}_i), \qquad \|\phi\|_{2,n}^2 = n^{-1}\sum_{i=1}^n \phi^2(\mathbf{X}_i).$$

In addition, if the functions $\phi, \varphi$ are $L^2$-integrable, define the theoretical inner product and its corresponding theoretical $L^2$ norm as

$$\langle \phi, \varphi \rangle_2 = E\{\phi(\mathbf{X}_i)\varphi(\mathbf{X}_i)\}, \qquad \|\phi\|_2^2 = E\{\phi^2(\mathbf{X}_i)\}.$$

The evaluation of spline estimator $\hat{m}(\mathbf{x})$ at the $n$ observations results in an $n$-dimensional vector, $\hat{m}(\mathbf{X}_1,\ldots,\mathbf{X}_n) = \{\hat{m}(\mathbf{X}_1),\ldots,\hat{m}(\mathbf{X}_n)\}^T$, which can be considered as the projection of $\mathbf{Y}$ on the space $G_n$ with respect to the empirical inner product $\langle \cdot, \cdot \rangle_{2,n}$. In general, for any $n$-dimensional vector $\mathbf{\Lambda} = \{\Lambda_1,\ldots,\Lambda_n\}^T$, we define $\mathbf{P}_n\mathbf{\Lambda}(\mathbf{x})$ as the spline function constructed from the projection of $\mathbf{\Lambda}$ on the inner product space $(G_n, \langle \cdot, \cdot \rangle_{2,n})$, that is,

$$\mathbf{P}_n\mathbf{\Lambda}(\mathbf{x}) = \hat{\lambda}_0 + \sum_{\alpha=1}^d \sum_{J=1}^N \hat{\lambda}_{J,\alpha} B_{J,\alpha}(x_\alpha),$$

with the coefficients $(\hat{\lambda}_0, \hat{\lambda}_{1,1}, \ldots, \hat{\lambda}_{N,d})$ given in (2.6). Next, the multivariate function $\mathbf{P}_n \mathbf{\Lambda}(\mathbf{x})$ is decomposed into the empirically centered additive components $\mathbf{P}_{n,\alpha} \mathbf{\Lambda}(x_\alpha)$, $\alpha = 1, \ldots, d$, and the constant component $\mathbf{P}_{n,c} \mathbf{\Lambda}$:

$$(3.1) \qquad \mathbf{P}_{n,\alpha} \mathbf{\Lambda}(x_\alpha) = \mathbf{P}^*_{n,\alpha} \mathbf{\Lambda}(x_\alpha) - n^{-1} \sum_{i=1}^{n} \mathbf{P}^*_{n,\alpha} \mathbf{\Lambda}(X_{i\alpha}),$$

$$(3.2) \qquad \mathbf{P}_{n,c} \mathbf{\Lambda} = \hat{\lambda}_0 + n^{-1} \sum_{\alpha=1}^{d} \sum_{i=1}^{n} \mathbf{P}^*_{n,\alpha} \mathbf{\Lambda}(X_{i\alpha}),$$

where $\mathbf{P}^*_{n,\alpha} \mathbf{\Lambda}(x_\alpha) = \sum_{J=1}^{N} \hat{\lambda}_{J,\alpha} B_{J,\alpha}(x_\alpha)$. With this new notation, we can rewrite the spline estimators $\hat{m}(\mathbf{x}), \hat{m}_\alpha(x_\alpha), \hat{m}_c$ defined in (2.5) and (2.7) as

$$\hat{m}(\mathbf{x}) = \mathbf{P}_n \mathbf{Y}(\mathbf{x}), \qquad \hat{m}_\alpha(x_\alpha) = \mathbf{P}_{n,\alpha} \mathbf{Y}(x_\alpha), \qquad \hat{m}_c = \mathbf{P}_{n,c} \mathbf{Y}.$$

Based on the relation $Y_i = m(\mathbf{X}_i) + \sigma(\mathbf{X}_i)\varepsilon_i$, one defines similarly the noiseless spline smoothers and the variance spline components,

$$(3.3) \qquad \begin{aligned} \tilde{m}(\mathbf{x}) &= \mathbf{P}_n\{m(\mathbf{X})\}(\mathbf{x}), \qquad \tilde{m}_\alpha(x_\alpha) = \mathbf{P}_{n,\alpha}\{m(\mathbf{X})\}(x_\alpha), \\ \tilde{m}_c &= \mathbf{P}_{n,c}\{m(\mathbf{X})\}, \end{aligned}$$

$$(3.4) \qquad \tilde{\varepsilon}(\mathbf{x}) = \mathbf{P}_n \mathbf{E}(\mathbf{x}), \qquad \tilde{\varepsilon}_\alpha(x_\alpha) = \mathbf{P}_{n,\alpha} \mathbf{E}(x_\alpha), \qquad \tilde{\varepsilon}_c = \mathbf{P}_{n,c} \mathbf{E},$$

where the noise vector $\mathbf{E} = \{\sigma(\mathbf{X}_i)\varepsilon_i\}_{i=1}^{n}$. Due to the linearity of the operators $\mathbf{P}_n$, $\mathbf{P}_{n,c}$, $\mathbf{P}_{n,\alpha}$, $\alpha = 1, \ldots, d$, one has the crucial decomposition

$$(3.5) \qquad \begin{aligned} \hat{m}(\mathbf{x}) &= \tilde{m}(\mathbf{x}) + \tilde{\varepsilon}(\mathbf{x}), \qquad \hat{m}_c = \tilde{m}_c + \tilde{\varepsilon}_c, \\ \hat{m}_\alpha(x_\alpha) &= \tilde{m}_\alpha(x_\alpha) + \tilde{\varepsilon}_\alpha(x_\alpha), \end{aligned}$$

for $\alpha = 1, \ldots, d$. As closer examination is needed later for $\tilde{\varepsilon}(\mathbf{x})$ and $\tilde{\varepsilon}_\alpha(x_\alpha)$, we define in addition $\tilde{\mathbf{a}} = \{\tilde{a}_0, \tilde{a}_{1,1}, \ldots, \tilde{a}_{N,d}\}^T$ as the minimizer of

$$(3.6) \qquad \sum_{i=1}^{n} \left\{ \sigma(\mathbf{X}_i)\varepsilon_i - a_0 - \sum_{\alpha=1}^{d} \sum_{J=1}^{N} a_{J,\alpha} B_{J,\alpha}(X_{i\alpha}) \right\}^2.$$

Then $\tilde{\varepsilon}(\mathbf{x}) = \tilde{\mathbf{a}}^T \mathbf{B}(\mathbf{x})$, where the vector $\mathbf{B}(\mathbf{x})$ and matrix $\mathbf{B}$ are defined as

$$(3.7) \quad \mathbf{B}(\mathbf{x}) = \{1, B_{1,1}(x_1), \ldots, B_{N,d}(x_d)\}^T, \qquad \mathbf{B} = \{\mathbf{B}(\mathbf{X}_1), \ldots, \mathbf{B}(\mathbf{X}_n)\}^T.$$

Thus $\tilde{\mathbf{a}} = (\mathbf{B}^T \mathbf{B})^{-1} \mathbf{B}^T \mathbf{E}$ is the solution of (3.6) and specifically $\tilde{\mathbf{a}}$ is equal to

$$(3.8) \qquad \begin{aligned} &\begin{Bmatrix} 1 & \mathbf{0}_{dN}^T \\ \mathbf{0}_{dN} & \langle B_{J,\alpha}, B_{J',\alpha'} \rangle_{2,n} \end{Bmatrix}_{\substack{1 \le \alpha, \alpha' \le d, \\ 1 \le J, J' \le N}}^{-1} \\ &\qquad \times \begin{Bmatrix} \dfrac{1}{n} \sum_{i=1}^{n} \sigma(\mathbf{X}_i)\varepsilon_i \\ \dfrac{1}{n} \sum_{i=1}^{n} B_{J,\alpha}(X_{i\alpha}) \sigma(\mathbf{X}_i)\varepsilon_i \end{Bmatrix}_{\substack{1 \le J \le N, \\ 1 \le \alpha \le d}}, \end{aligned}$$

where $\mathbf{0}_p$ is a $p$-vector with all elements 0.

Our main objective is to study the difference between the smoothed backfitted estimator $\hat{m}_1^*(x_1)$ and the smoothed "oracle" estimator $\tilde{m}_1^*(x_1)$, both given in (2.9). From now on, we assume without loss of generality that $d=2$ for notational brevity. Making use of the definition of $\hat{c}$ and the signal and noise decomposition (3.5), the difference $\tilde{m}_1^*(x_1) - \hat{m}_1^*(x_1) - \hat{c} + c$ can be treated as the sum of two terms,

$$
(3.9) \quad \frac{1/n \sum_{i=1}^n K_h(X_{i1}-x_1)\{\hat{m}_2(X_{i2}) - m_2(X_{i2})\}}{1/n \sum_{i=1}^n K_h(X_{i1}-x_1)} = \frac{\Psi_b(x_1) + \Psi_v(x_1)}{1/n \sum_{i=1}^n K_h(X_{i1}-x_1)},
$$

where

$$
(3.10) \quad \Psi_b(x_1) = \frac{1}{n} \sum_{i=1}^n K_h(X_{i1}-x_1)\{\tilde{m}_2(X_{i2}) - m_2(X_{i2})\},
$$

$$
(3.11) \quad \Psi_v(x_1) = \frac{1}{n} \sum_{i=1}^n K_h(X_{i1}-x_1)\tilde{\varepsilon}_2(X_{i2}).
$$

The term $\Psi_b(x_1)$ is induced by the bias term $\tilde{m}_2(X_{i2}) - m_2(X_{i2})$, while $\Psi_v(x_1)$ is related to the variance term $\tilde{\varepsilon}_2(X_{i2})$. Both of these two terms have order $o_p(n^{-2/5})$ by Propositions 4.1 and 5.1 in the next two sections. Standard theory of kernel density estimation ensures that the denominator in (3.9), $n^{-1} \sum_{i=1}^n K_h(X_{i1}-x_1)$, has a positive lower bound for $x_1 \in [0,1]$. The additional nuisance term $\hat{c} - c$ is clearly of order $O_p(n^{-1/2})$ and thus $o_p(n^{-2/5})$, which needs no further arguments for the proofs. Theorem 2.1 then follows from Propositions 4.1 and 5.1.

**4. Bias reduction for $\Psi_b(x_1)$.** In this section we show that the bias term $\Psi_b(x_1)$ of (3.10) is uniformly of order $o_p(n^{-2/5})$ for $x_1 \in [0,1]$.

PROPOSITION 4.1. *Under assumptions* (A1), (A2) *and* (A4)–(A6),

$$
\sup_{x_1 \in [0,1]} |\Psi_b(x_1)| = O_p(n^{-1/2} + H) = o_p(n^{-2/5}).
$$

LEMMA 4.1. *Under assumption* (A1), *there exist functions* $g_1$, $g_2 \in G$, *such that*

$$
\left\| \tilde{m} - g + \sum_{\alpha=1}^2 \langle 1, g_\alpha(X_\alpha) \rangle_{2,n} \right\|_{2,n} = O_p(n^{-1/2} + H),
$$

*where* $g(\mathbf{x}) = c + \sum_{\alpha=1}^2 g_\alpha(x_\alpha)$ *and* $\tilde{m}$ *is defined in* (3.3).

PROOF. According to the result on page 149 of [3], there is a constant $C_\infty > 0$ such that for the function $g_\alpha \in G$, $\|g_\alpha - m_\alpha\|_\infty \leq C_\infty H$, $\alpha = 1, 2$. Thus $\|g - m\|_\infty \leq \sum_{\alpha=1}^{2} \|g_\alpha - m_\alpha\|_\infty \leq 2C_\infty H$ and $\|\tilde{m} - m\|_{2,n} \leq \|g - m\|_{2,n} \leq 2C_\infty H$. Noting that $\|\tilde{m} - g\|_{2,n} \leq \|\tilde{m} - m\|_{2,n} + \|g - m\|_{2,n} \leq 4C_\infty H$, one has

$$|\langle g_\alpha(X_\alpha), 1\rangle_{2,n}| \leq |\langle 1, g_\alpha(X_\alpha)\rangle_{2,n} - \langle 1, m_\alpha(X_\alpha)\rangle_{2,n}| + |\langle 1, m_\alpha(X_\alpha)\rangle_{2,n}| \quad (4.1)$$
$$\leq C_\infty H + O_p(n^{-1/2}).$$

Therefore

$$\left\|\tilde{m} - g + \sum_{\alpha=1}^{2} \langle 1, g_\alpha(X_\alpha)\rangle_{2,n}\right\|_{2,n} \leq \|\tilde{m} - g\|_{2,n} + \sum_{\alpha=1}^{2} |\langle 1, g_\alpha(X_\alpha)\rangle_{2,n}|$$
$$\leq 6C_\infty H + O_p(n^{-1/2}) = O_p(n^{-1/2} + H).$$

□

PROOF OF PROPOSITION 4.1. Denote

$$R_1 = \sup_{x_1 \in [0,1]} \left|\frac{\sum_{i=1}^{n} K_h(X_{i1} - x_1)\{g_2(X_{i2}) - m_2(X_{i2})\}}{\sum_{i=1}^{n} K_h(X_{i1} - x_1)}\right|,$$

$$R_2 = \sup_{x_1 \in [0,1]} \left|\frac{\sum_{i=1}^{n} K_h(X_{i1} - x_1)\{\tilde{m}_2(X_{i2}) - g_2(X_{i2}) + \langle 1, g_2(X_2)\rangle_{2,n}\}}{\sum_{i=1}^{n} K_h(X_{i1} - x_1)}\right|;$$

then $\sup_{x_1 \in [0,1]} |\Psi_b(x_1)| \leq |\langle 1, g_2(X_2)\rangle_{2,n}| + R_1 + R_2$. For $R_1$, using the result on page 149 of [3], one has $R_1 \leq C_\infty H$. To deal with $R_2$, let $B_{J,2}^*(x_\alpha) = B_{J,2}(x_\alpha) - \langle 1, B_{J,2}(X_\alpha)\rangle_{2,n}$, for $J = 1, \ldots, N$, $\alpha = 1, 2$; then one can write

$$\tilde{m}(\mathbf{x}) - g(\mathbf{x}) + \sum_{\alpha=1}^{2} \langle 1, g_\alpha(X_\alpha)\rangle_{2,n} = \tilde{a}^* + \sum_{\alpha=1}^{2} \sum_{J=1}^{N} \tilde{a}_{J,\alpha}^* B_{J,\alpha}^*(x_\alpha).$$

Thus, $n^{-1} \sum_{i=1}^{n} K_h(X_{i1} - x_1)\{\tilde{m}_2(X_{i2}) - g_2(X_{i2}) + \langle 1, g_2(X_2)\rangle_{2,n}\}$ can be rewritten as $n^{-1} \sum_{i=1}^{n} K_h(X_{i1} - x_1) \sum_{J=1}^{N} \tilde{a}_{J,2}^* B_{J,2}^*(X_{i2})$, bounded by

$$\sum_{J=1}^{N} |\tilde{a}_{J,2}^*| \sup_{1 \leq J \leq N} \left|n^{-1} \sum_{i=1}^{n} K_h(X_{i1} - x_1) B_{J,2}^*(X_{i2})\right|$$

$$\leq \sum_{J=1}^{N} |\tilde{a}_{J,2}^*| \left\{\sup_{1 \leq J \leq N} \left|n^{-1} \sum_{i=1}^{n} \omega_J(\mathbf{X}_i, x_1)\right| + A_{n,1} \left|n^{-1} \sum_{i=1}^{n} K_h(X_{i1} - x_1)\right|\right\},$$

where $A_{n,1} = \sup_{J,\alpha} |\langle 1, B_{J,\alpha}\rangle_{2,n} - \langle 1, B_{J,\alpha}\rangle_2| = O_p(n^{-1/2} \log n)$ as in (A.12) and $\omega_J(\mathbf{X}_i, x_1)$ is in (5.5) with mean $\mu_{\omega_J}(x_1)$. By Lemma A.3

$$\sup_{x_1 \in [0,1]} \sup_{1 \leq J \leq N} \left|\frac{1}{n}\sum_{i=1}^{n} \omega_J(\mathbf{X}_i, x_1)\right|$$

$$\leq \sup_{x_1 \in [0,1]} \sup_{1 \leq J \leq N} \left| \frac{1}{n} \sum_{i=1}^{n} \omega_J(\mathbf{X}_i, x_1) - \mu_{\omega_J}(x_1) \right|$$

$$+ \sup_{x_1 \in [0,1]} \sup_{1 \leq J \leq N} |\mu_{\omega_J}(x_1)|$$

$$= O_p(\log n/\sqrt{nh}) + O_p(H^{1/2}) = O_p(H^{1/2}).$$

Therefore, one has

$$\sup_{x_1 \in [0,1]} \left| n^{-1} \sum_{i=1}^{n} K_h(X_{i1} - x_1) \{ \tilde{m}_2(X_{i2}) - g_2(X_{i2}) + \langle 1, g_2(X_2) \rangle_{2,n} \} \right|$$

$$\leq \left\{ N \sum_{J=1}^{N} (\tilde{a}_{J,2}^*)^2 \right\}^{1/2} \left\{ O_p(H^{1/2}) + O_p\left(\frac{\log n}{\sqrt{n}}\right) \right\} = O_p\left( \left\{ \sum_{J=1}^{N} (\tilde{a}_{J,2}^*)^2 \right\}^{1/2} \right)$$

$$= O_p\left( \left\| \tilde{m} - g + \sum_{\alpha=1}^{2} \langle 1, g_\alpha(X_\alpha) \rangle_{2,n} \right\|_2 \right)$$

$$= O_p\left( \left\| \tilde{m} - g + \sum_{\alpha=1}^{2} \langle 1, g_\alpha(X_\alpha) \rangle_{2,n} \right\|_{2,n} \right),$$

where the last step follows from Lemma A.8. Thus, by Lemma 4.1,

(4.2) $$R_2 = O_p(n^{-1/2} + H).$$

Combining (4.1) and (4.2), one establishes Proposition 4.1. □

**5. Variance reduction for $\Psi_v(x_1)$.** In this section we will see that the term $\Psi_v(x_1)$ given in (3.11) is uniformly of order $o_p(n^{-2/5})$. This is the most challenging part to be proved, mostly done in the Appendix. Define an auxiliary entity

(5.1) $$\tilde{\varepsilon}_2^* = \sum_{J=1}^{N} \tilde{a}_{J,2} B_{J,2}(x_2),$$

where $\tilde{a}_{J,2}$ is given in (3.8). Definitions (3.1) and (3.2) imply that $\tilde{\varepsilon}_2(x_2)$ defined in (3.4) is simply the empirical centering of $\tilde{\varepsilon}_2^*(x_2)$, that is,

(5.2) $$\tilde{\varepsilon}_2(x_2) \equiv \tilde{\varepsilon}_2^*(x_2) - n^{-1} \sum_{i=1}^{n} \tilde{\varepsilon}_2^*(X_{i2}).$$

PROPOSITION 5.1. *Under assumptions* (A2)–(A6), *one has*

$$\sup_{x_1 \in [0,1]} |\Psi_v(x_1)| = O_p(H) = o_p(n^{-2/5}).$$

According to (5.2), we can write $\Psi_v(x_1) = \Psi_v^{(2)}(x_1) - \Psi_v^{(1)}(x_1)$, where

$$\Psi_v^{(1)}(x_1) = n^{-1} \sum_{l=1}^n K_h(X_{l1} - x_1) \cdot n^{-1} \sum_{i=1}^n \tilde{\varepsilon}_2^*(X_{i2}), \tag{5.3}$$

$$\Psi_v^{(2)}(x_1) = n^{-1} \sum_{l=1}^n K_h(X_{l1} - x_1) \tilde{\varepsilon}_2^*(X_{l2}), \tag{5.4}$$

in which $\tilde{\varepsilon}_2^*(X_{i2})$ is given in (5.1). Further one denotes

$$\omega_J(\mathbf{X}_l, x_1) = K_h(X_{l1} - x_1) B_{J,2}(X_{l2}), \qquad \mu_{\omega_J}(x_1) = E\omega_J(\mathbf{X}_l, x_1). \tag{5.5}$$

By (3.8) and (5.1), $\Psi_v^{(2)}(x_1)$ can be rewritten as

$$\Psi_v^{(2)}(x_1) = n^{-1} \sum_{l=1}^n \sum_{J=1}^N \tilde{a}_{J,2} \omega_J(\mathbf{X}_l, x_1). \tag{5.6}$$

The uniform order of $\Psi_v^{(1)}(x_1)$ and $\Psi_v^{(2)}(x_1)$ is given in the next two lemmas.

LEMMA 5.1. *Under assumptions* (A2)–(A6), $\Psi_v^{(1)}(x_1)$ *in (5.3) satisfies*
$$\sup_{x_1 \in [0,1]} |\Psi_v^{(1)}(x_1)| = O_p\{N(\log n)^2/n\}.$$

PROOF. Based on (5.1),
$$n^{-1} \sum_{i=1}^n \tilde{\varepsilon}_2^*(X_{i2}) \le \left|\sum_{J=1}^N \tilde{a}_{J,2}\right| \cdot \sup_{1 \le J \le N} \left|\frac{1}{n} \sum_{i=1}^n B_{J,2}(X_{i2})\right|.$$

Lemma A.6 implies that
$$\left|\sum_{J=1}^N \tilde{a}_{J,2}\right| \le \left\{N \sum_{J=1}^N \tilde{a}_{J,2}^2\right\}^{1/2} \le \{N\tilde{\mathbf{a}}^T \tilde{\mathbf{a}}\}^{1/2} = O_p(Nn^{-1/2} \log n).$$

By (A.12), $\sup_{1 \le J \le N} |n^{-1} \sum_{i=1}^n B_{J,2}(X_{i2})| \le A_{n,1} = O_p(n^{-1/2} \log n)$, so

$$\frac{1}{n} \sum_{i=1}^n \tilde{\varepsilon}_2^*(X_{i2}) = O_p\{N(\log n)^2/n\}. \tag{5.7}$$

By assumption (A5) on the kernel function $K$, standard theory on kernel density estimation entails that $\sup_{x_1 \in [0,1]} |n^{-1} \sum_{l=1}^n K_h(X_{l1} - x_1)| = O_p(1)$. Thus with (5.7) the lemma follows immediately. □

LEMMA 5.2. *Under assumptions* (A2)–(A6), $\Psi_v^{(2)}(x_1)$ *in (5.4) satisfies*
$$\sup_{x_1 \in [0,1]} |\Psi_v^{(2)}(x_1)| = O_p(H).$$

Lemma 5.2 follows from Lemmas A.10 and A.11. Proposition 5.1 follows from Lemmas 5.1 and 5.2.

**6. Simulation example.** In this section we carry out two simulation experiments to illustrate the finite-sample behavior of our SPBK estimators. The programming codes are available in both R 2.2.1 and XploRe. For information on XploRe, see [8] or visit www.xplore-stat.de.

The number of interior knots $N$ for the spline estimation as in (2.6) will be determined by the sample size $n$ and a tuning constant $c$. To be precise,

$$N = \min([cn^{2/5} \log n] + 1, [(n/2 - 1)d^{-1}]),$$

in which $[a]$ denotes the integer part of $a$. In our simulation study, we have used $c = 0.5, 1.0$. As seen in Table 1, the choice of $c$ makes little difference, so we always recommend to use $c = 0.5$ to save computation for massive data sets. The additional constraint that $N \leq (n/2 - 1)d^{-1}$ ensures that the number of terms in the linear least squares problem (2.6), $1 + dN$, is no greater than $n/2$, which is necessary when the sample size $n$ is moderate and the dimension $d$ is high.

We have obtained for comparison both the SPBK estimator $\hat{m}_\alpha^*(x_\alpha)$ and the "oracle" estimator $\tilde{m}_\alpha^*(x_\alpha)$ by Nadaraya–Watson regression estimation using a quartic kernel and the rule-of-thumb bandwidth.

We consider first the accuracy of the estimation, measured in terms of mean average squared error. Then to see that the SPBK estimator $\hat{m}_\alpha^*(x_\alpha)$ is as efficient as the "oracle smoother" $\tilde{m}_\alpha^*(x_\alpha)$, we define the empirical relative efficiency of $\hat{m}_\alpha^*(x_\alpha)$ with respect to $\tilde{m}_\alpha^*(x_\alpha)$ as

$$(6.1) \quad \text{eff}_\alpha = \left[\frac{\sum_{i=1}^n \{\tilde{m}_\alpha^*(X_{i\alpha}) - m_\alpha(X_{i\alpha})\}^2}{\sum_{i=1}^n \{\hat{m}_\alpha^*(X_{i\alpha}) - m_\alpha(X_{i\alpha})\}^2}\right]^{1/2}.$$

Theorem 2.1 indicates that the $\text{eff}_\alpha$ should be close to 1 for all $\alpha = 1, \ldots, d$. Figure 2 provides the kernel density estimations of the above empirical efficiencies to observe the convergence.

EXAMPLE 6.1. A time series $\{Y_t\}_{t=-1999}^{n+3}$ is generated according to the nonlinear additive autoregression model with sine functions given in [2],

$$Y_t = 1.5 \sin\left(\frac{\pi}{2} Y_{t-2}\right) - 1.0 \sin\left(\frac{\pi}{2} Y_{t-3}\right) + \sigma_0 \varepsilon_t, \qquad \sigma_0 = 0.5, 1.0,$$

where $\{\varepsilon_t\}_{t=-1996}^{n+3}$ are i.i.d. standard normal errors. Let $\mathbf{X}_t^T = \{Y_{t-1}, Y_{t-2}, Y_{t-3}\}$. Theorem 3 on page 91 of [4] establishes that $\{Y_t, \mathbf{X}_t^T\}_{t=-1996}^{n+3}$ is geometrically ergodic. The first 2000 observations are discarded to make $\{Y_t\}_{t=1}^{n+3}$ behave like a geometrically $\alpha$-mixing and strictly stationary time series. The multivariate datum $\{Y_t, \mathbf{X}_t^T\}_{t=4}^{n+3}$ then satisfies assumptions (A1) to (A6) except

ADDITIVE AUTOREGRESSION MODEL                     15TABLE 1
*Report of Example 6.1*

| $\sigma_0$ | $n$ | $c$ | Component #1 | | Component #2 | | Component #3 | |
|---|---|---|---|---|---|---|---|---|
| | | | 1st stage | 2nd stage | 1st stage | 2nd stage | 1st stage | 2nd stage |
| 0.5 | 100 | 0.5 | 0.1231 | 0.0461 | 0.1476 | 0.0645 | 0.1254 | 0.0681 |
| | | 1.0 | 0.1278 | 0.0520 | 0.1404 | 0.0690 | 0.1318 | 0.0726 |
| | 200 | 0.5 | 0.0539 | 0.0125 | 0.0616 | 0.0275 | 0.0577 | 0.0252 |
| | | 1.0 | 0.0841 | 0.0144 | 0.0839 | 0.0290 | 0.0848 | 0.0285 |
| | 500 | 0.5 | 0.0263 | 0.0031 | 0.0306 | 0.0107 | 0.0278 | 0.0102 |
| | | 1.0 | 0.0595 | 0.0044 | 0.0578 | 0.0115 | 0.0605 | 0.0119 |
| | 1000 | 0.5 | 0.0169 | 0.0015 | 0.0210 | 0.0053 | 0.0178 | 0.0054 |
| | | 1.0 | 0.0364 | 0.0018 | 0.0367 | 0.0054 | 0.0375 | 0.0059 |
| 1.0 | 100 | 0.5 | 0.3008 | 0.0587 | 0.3298 | 0.1427 | 0.3236 | 0.1393 |
| | | 1.0 | 0.3088 | 0.0586 | 0.3369 | 0.1364 | 0.3062 | 0.1316 |
| | 200 | 0.5 | 0.1742 | 0.0256 | 0.1783 | 0.0802 | 0.1892 | 0.0701 |
| | | 1.0 | 0.2899 | 0.0328 | 0.2830 | 0.0824 | 0.3043 | 0.0721 |
| | 500 | 0.5 | 0.0924 | 0.0065 | 0.1124 | 0.0421 | 0.1004 | 0.0345 |
| | | 1.0 | 0.2299 | 0.0078 | 0.2305 | 0.0458 | 0.2314 | 0.0362 |
| | 1000 | 0.5 | 0.0616 | 0.0033 | 0.0637 | 0.0270 | 0.0646 | 0.0224 |
| | | 1.0 | 0.1460 | 0.0034 | 0.1433 | 0.0275 | 0.1429 | 0.0219 |

Monte Carlo average squared errors (ASE) based on 100 replications.

that instead of being $[0,1]$, the range of $Y_{t-\alpha}$, $\alpha = 1, 2, 3$, needs to be recalibrated. Since we have no exact knowledge of the distribution of the $Y_t$, we have generated many realizations of size 50,000 from which we found that more than 95% of the observations fall in $[-2.58, 2.58]$ ($[-3.14, 3.14]$) with $\sigma_0 = 0.5$ ($\sigma_0 = 1$). We will estimate the functions $\{m_\alpha(x_\alpha)\}_{\alpha=1}^3$ for $x_\alpha \in [-2.58, 2.58]$ ($[-3.14, 3.14]$) with $\sigma_0 = 0.5$ ($\sigma_0 = 1.0$), where

$$m_1(x_1) \equiv 0, \qquad m_2(x_2) \equiv 1.5 \sin\left(\frac{\pi}{2} x_2\right) - E\left[1.5 \sin\left(\frac{\pi}{2} Y_t\right)\right],$$

(6.2)

$$m_3(x_3) \equiv -1.0 \sin\left(\frac{\pi}{2} x_3\right) - E\left[-1.0 \sin\left(\frac{\pi}{2} Y_t\right)\right].$$

We choose the sample size $n$ to be 100, 200, 500 and 1000. Table 1 lists the average squared error (ASE) of the SPBK estimators and the constant spline pilot estimators from 100 Monte Carlo replications. As expected, increases in sample size reduce ASE for both estimators and across all combinations of $c$ values and noise levels. Table 1 also shows that our SPBK estimators improve upon the spline pilot estimators immensely regardless of noise level and sample size, which implies that our second Nadaraya–Watson smoothing step is not redundant.

To have some impression of the actual function estimates, at noise level $\sigma_0 = 0.5$ with sample size $n = 200, 500$, we have plotted the oracle estimator (thin dotted lines), SPBK estimator $\hat{m}_\alpha^*$ (thin solid lines) and their 95% pointwise confidence intervals (upper and lower dashed curves) for the true functions $m_\alpha$ (thick solid lines) in Figure 1. The visual impression of the SPBK estimators is rather satisfactory and their performance improves with increasing $n$.

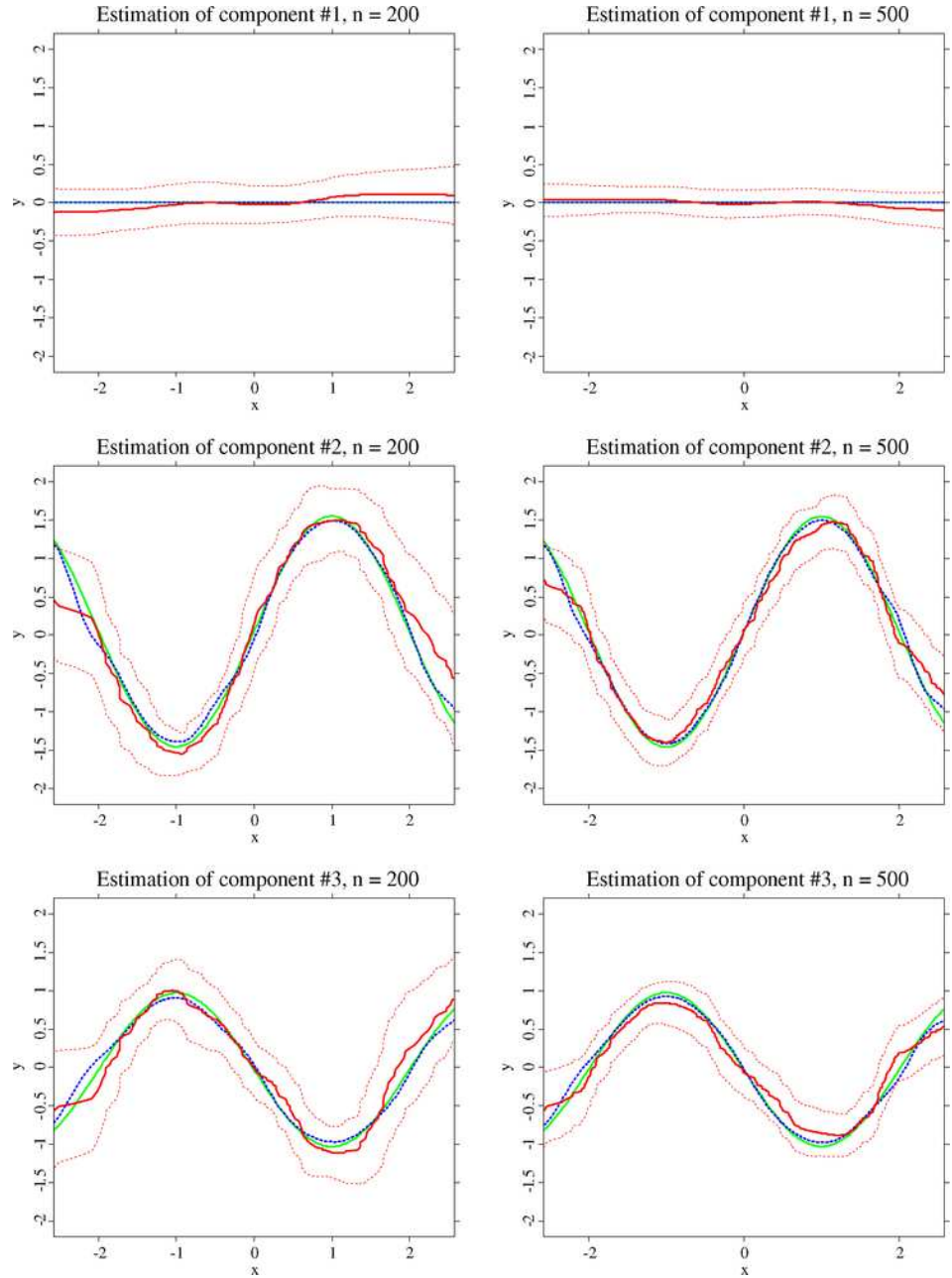

FIG. 1. *Plots of the oracle estimator (dotted blue curve), SPBK estimator (solid red curve) and the 95% pointwise confidence intervals constructed by (2.13) (upper and lower dashed red curves) of the function components $m_\alpha(x_\alpha)$ in (6.2), $\alpha = 1, 2, 3$ (solid green curve).*

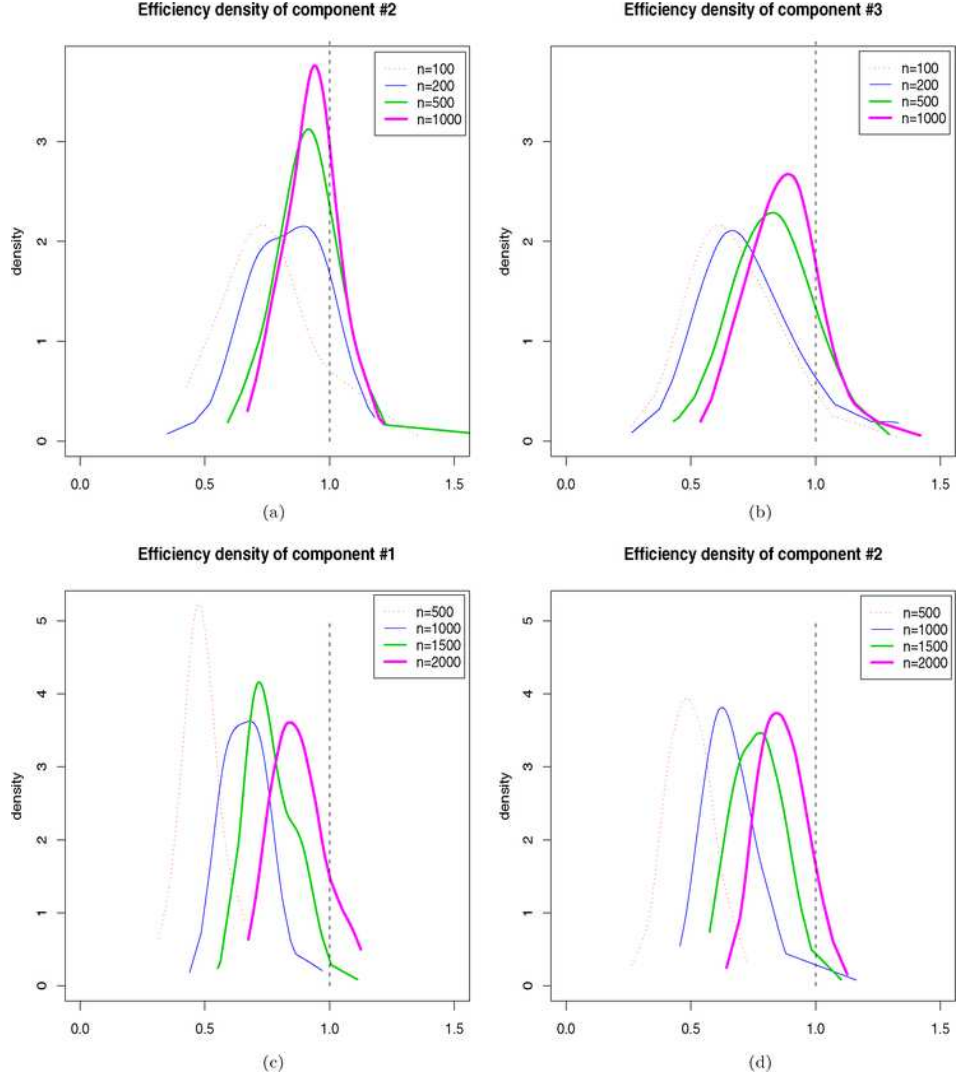

Fig. 2. *Kernel density plots of the* 100 *empirical efficiencies of* $\hat{m}_\alpha^*(x_\alpha)$ *to* $\tilde{m}_\alpha^*(x_\alpha)$, *computed according to (6.1):* (a) *Example 6.1* ($\alpha = 2, d = 3$); (b) *Example 6.1* ($\alpha = 3, d = 3$); (c) *Example 6.2* ($\alpha = 1, d = 30$); (d) *Example 6.2* ($\alpha = 2, d = 30$).

To see the convergence, Figure 2(a) and (b) plots the kernel density estimation of the 100 empirical efficiencies for $\alpha = 2, 3$ and sample sizes $n = 100, 200, 500$ and $1000$ at the noise level $\sigma_0 = 0.5$. The vertical line at efficiency $= 1$ is the standard line for the comparison of $\hat{m}_\alpha^*(x_\alpha)$ and $\tilde{m}_\alpha^*(x_\alpha)$. One can clearly see that the center of the density plots is going toward the standard line 1.0 with narrower spread when sample size $n$ is increasing, which is confirmative to the result of Theorem 2.1.

EXAMPLE 6.2. Consider the nonlinear additive heteroscedastic model

$$Y_t = \sum_{\alpha=1}^{d} \sin\left(\frac{\pi}{2.5} X_{t-\alpha}\right) + \sigma(\mathbf{X})\varepsilon_t, \qquad \varepsilon_t \stackrel{\text{i.i.d.}}{\sim} N(0,1),$$

in which $\mathbf{X}_t^T = \{X_{t-1}, \ldots, X_{t-d}\}$ is a sequence of i.i.d. standard normal random variables truncated by $[-2.5, 2.5]$ and

$$\sigma(\mathbf{X}) = \sigma_0 \frac{\sqrt{d}}{2} \cdot \frac{5 - \exp(\sum_{\alpha=1}^{d} |X_{t-\alpha}|/d)}{5 + \exp(\sum_{\alpha=1}^{d} |X_{t-\alpha}|/d)}, \qquad \sigma_0 = 0.1.$$

By this choice of $\sigma(\mathbf{X})$, we ensure that our design is heteroscedastic, and the variance is roughly proportional to dimension $d$, which is intended to mimic the case when independent copies of the same kind of univariate regression problem are simply added together.

For $d = 30$, we have run 100 replications for sample size $n = 500, 1000, 1500$ and $2000$. The kernel density estimation of the 100 empirical efficiencies for $\alpha = 1, 2$ is graphically represented respectively in (c) and (d) of Figure 2. Again one sees that with increasing $n$, the efficiency distribution converges to 1.

Lastly, we provide the computing time of Example 6.1 from 100 replications on an ordinary PC with Intel Pentium IV 1.86 GHz processor and 1.0 GB RAM. The average time run by XploRe to generate one sample of size $n$ and compute the SPBK estimator and marginal integration estimator (MIE) is reported in Table 2. The MIEs have been obtained by directly recalling the "intest" in XploRe. As expected, the computing time for MIE is extremely sensitive to sample size due to the fact that it requires $n^2$ least squares in two steps. In contrast, at least for large sample data, our proposed SPBK is thousands of times faster than MIE. Thus our SPBK estimation is feasible and appealing to deal with massive data sets.

## APPENDIX

Throughout this section, $a_n \gg b_n$ means $\lim_{n \to \infty} b_n/a_n = 0$ and $a_n \sim b_n$ means $\lim_{n \to \infty} b_n/a_n = c$, where $c$ is some constant.

ADDITIVE AUTOREGRESSION MODEL    19TABLE 2
*The computing time of Example 6.1 (in seconds)*

| **Method** | $n = 100$ | $n = 200$ | $n = 400$ | $n = 1000$ |
|---|---|---|---|---|
| MIE | 10 | 76 | 628 | 10728 |
| SPBK | 0.7 | 0.9 | 1.2 | 4.5 |

**A.1. Preliminaries.** We first give the Bernstein inequality for a geometric $\alpha$-mixing sequence, which plays an important role through our proof.

LEMMA A.1 (Theorem 1.4, page 31 of [1]). *Let $\{\xi_t, t \in \mathbb{Z}\}$ be a zero-mean real-valued $\alpha$-mixing process, $S_n = \sum_{i=1}^n \xi_i$. Suppose that there exists $c > 0$ such that for $i = 1, \ldots, n$, $k = 3, 4, \ldots, E|\xi_i|^k \le c^{k-2} k! E \xi_i^2 < +\infty$; then for each $n > 1$, integer $q \in [1, n/2]$, each $\varepsilon > 0$ and $k \ge 3$,*

$$P(|S_n| \ge n\varepsilon) \le a_1 \exp\left(-\frac{q\varepsilon^2}{25 m_2^2 + 5 c \varepsilon}\right) + a_2(k) \alpha\left(\left[\frac{n}{q+1}\right]\right)^{2k/(2k+1)},$$

*where $\alpha(\cdot)$ is the $\alpha$-mixing coefficient defined in (2.10) and*

$$a_1 = 2\frac{n}{q} + 2\left(1 + \frac{\varepsilon^2}{25 m_2^2 + 5 c \varepsilon}\right), \qquad a_2(k) = 11 n \left(1 + \frac{5 m_k^{2k/(2k+1)}}{\varepsilon}\right),$$

*with $m_r = \max_{1 \le i \le n} \|\xi_i\|_r$, $r \ge 2$.*

LEMMA A.2. *Under assumptions* (A4) *and* (A6), *one has:*

(i) *There exist constants $C_0(f)$ and $C_1(f)$ depending on the marginal densities $f_\alpha(x_\alpha)$, $\alpha = 1, 2$, such that $C_0(f) H \le \|b_{J,\alpha}\|_2^2 \le C_1(f) H$, where $b_{J,\alpha}$ is given in (2.2).*

(ii) *For any $\alpha = 1, 2$, $|J' - J| \le 1$, $E\{B_{J,\alpha}(X_{i\alpha}) B_{J',\alpha}(X_{i\alpha})\} \sim 1$; in addition*

$$E|B_{J,\alpha}(X_{i\alpha}) B_{J',\alpha}(X_{i\alpha})|^k \sim H^{1-k}, \qquad k \ge 1,$$

*where $B_{J,\alpha}$ and $B_{J',\alpha}$ are defined in (2.3).*

We refer the proof of the above lemma to Lemma A.2 in [26].

LEMMA A.3. *Under assumptions* (A4)–(A6), *for $\mu_{\omega_J}(x_1)$ given in (5.5),*

$$\sup_{x_1 \in [0,1]} \sup_{1 \le J \le N} |\mu_{\omega_J}(x_1)| = O(H^{1/2}).$$

PROOF. Denote the theoretical norm of $I_{J,\alpha}$ in (2.1) for $\alpha = 1, 2$, $J = 1, \ldots, N+1$,

$$(A.1) \qquad c_{J,\alpha} = \|I_{J,\alpha}\|_2^2 = \int I_{J,\alpha}^2(x_\alpha) f_\alpha(x_\alpha) \, dx_\alpha.$$

By definition, $|\mu_{\omega_J}(x_1)| = |E\{K_h(X_{l1} - x_1) B_{J,2}(X_{l2})\}|$ is bounded by

$$\iint K_h(u_1 - x_1) |B_{J,2}(u_2)| f(u_1, u_2) \, du_1 \, du_2$$

$$= (\|b_{J,2}\|_2)^{-1} \Big\{ \iint K(v_1) I_{J+1,2}(u_2) f(hv_1 + x_1, u_2) \, dv_1 \, du_2$$

$$+ \left(\frac{c_{J+1,2}}{c_{J,2}}\right)^{1/2} \iint K(v_1) I_{J,2}(u_2) f(hv_1 + x_1, u_2) \, dv_1 \, du_2 \Big\}.$$

The boundedness of the joint density $f$ and the Lipschitz continuity of the kernel $K$ will then imply that

$$\sup_{x_1 \in [0,1]} \sup_{1 \leq J \leq N} \iint K(v_1) I_{J,2}(u_2) f(hv_1 + x_1, u_2) \, dv_1 \, du_2 \leq C_K C_f H.$$

The proof of the lemma is then completed by (i) of Lemma A.2. □

LEMMA A.4. *Under assumptions* (A2) *and* (A4)–(A6), *one has*

$$(A.2) \quad \sup_{x_1 \in [0,1]} \sup_{1 \leq J \leq N} \left| n^{-1} \sum_{l=1}^n \{\omega_J(\mathbf{X}_l, x_1) - \mu_{\omega_J}(x_1)\} \right| = O_p(\log n / \sqrt{nh}),$$

$$(A.3) \qquad \sup_{x_1 \in [0,1]} \sup_{1 \leq J \leq N} \left| n^{-1} \sum_{l=1}^n \omega_J(\mathbf{X}_l, x_1) \right| = O_p(H^{1/2}),$$

*where $\omega_J(\mathbf{X}_l, x_1)$ and $\mu_{\omega_J}(x_1)$ are given in* (5.5).

PROOF. For simplicity, denote $\omega_J^*(\mathbf{X}_l, x_1) = \omega_J(\mathbf{X}_l, x_1) - \mu_{\omega_J}(x_1)$. Then

$$E\{\omega_J^*(\mathbf{X}_l, x_1)\}^2 = E\omega_J^2(\mathbf{X}_l, x_1) - \mu_{\omega_J}^2(x_1),$$

while $E\omega_J^2(\mathbf{X}_l, x_1)$ is equal to

$$h^{-1} \|b_{J,2}\|_2^{-2} \iint K^2(v_1) \Big\{ I_{J+1,2}(u_2) + \frac{c_{J+1,2}}{c_{J,2}} I_{J,2}(u_2) \Big\}$$

$$\times f(hv_1 + x_1, u_2) \, dv_1 \, du_2,$$

where $c_{J,\alpha}$ is given in (A.1). So $E\omega_J^2(\mathbf{X}_l, x_1) \sim h^{-1}$ and $E\omega_J^2(\mathbf{X}_l, x_1) \gg \mu_{\omega_J}^2(x_1)$. Hence for $n$ sufficiently large, $E\{\omega_J^*(\mathbf{X}_l, x_1)\}^2 = E\omega_J^2(\mathbf{X}_l, x_1) -$ $\mu^2_{\omega_J}(x_1) \geq c^* h^{-1}$, for some positive constant $c^*$. When $r \geq 3$, the $r$th moment $E|\omega_J(\mathbf{X}_l, x_1)|^r$ is

$$\frac{1}{\|b_{J,2}\|_2^r} \iint K_h^r(u_1 - x_1)\left\{I_{J+1,2}(u_2) + \left(\frac{c_{J+1,2}}{c_{J,2}}\right)^r I_{J,2}(u_2)\right\} f(u_1, u_2)\, du_1\, du_2.$$

It is clear that $E|\omega_J(\mathbf{X}_l, x_1)|^r \sim h^{(1-r)} H^{1-r/2}$. According to Lemma A.3, one has $|E\omega_J(\mathbf{X}_l, x_1)|^r \leq CH^{r/2}$, thus $E|\omega_J(\mathbf{X}_l, x_1)|^r \gg |\mu_{\omega_J}(x_1)|^r$. In addition

$$E|\omega_J^*(\mathbf{X}_l, x_1)|^r \leq \left\{\frac{c}{hH^{1/2}}\right\}^{(r-2)} r! E|\omega_J^*(\mathbf{X}_l, x_1)|^2,$$

so there exists $c_* = ch^{-1}H^{-1/2}$ such that $E|\omega_J^*(\mathbf{X}_l, x_1)|^r \leq c_*^{r-2} r! E|\omega_J^*(\mathbf{X}_l, x_1)|^2$, which implies that $\{\omega_J^*(\mathbf{X}_l, x_1)\}_{l=1}^n$ satisfies Cramér's condition. By Bernstein's inequality, for $r = 3$

$$P\left\{\left|\frac{1}{n}\sum_{l=1}^n \omega_J^*(\mathbf{X}_l, x_1)\right| \geq \rho_n\right\} \leq a_1 \exp\left(-\frac{q\rho_n^2}{25m_2^2 + 5c_*\rho_n}\right) + a_2(3)\alpha\left(\left[\frac{n}{q+1}\right]\right)^{6/7}$$

with $m_2^2 \sim h^{-1}$, $m_3 = \max_{1 \leq i \leq n} \|\omega_J^*(\mathbf{X}_l, x_1)\|_3 \leq \{C_0(2h^{-1})^2\}^{1/3}$ and

$$\rho_n = \rho \frac{\log n}{\sqrt{nh}}, \qquad a_1 = 2\frac{n}{q} + 2\left(1 + \frac{\rho_n^2}{25m_2^2 + 5c_*\rho_n}\right),$$

$$a_2(3) = 11n\left(1 + \frac{5m_3^{6/7}}{\rho_n}\right).$$

Observe that $5c_*\rho_n = o(1)$; by taking $q$ such that $[\frac{n}{q+1}] \geqslant c_0 \log n$, $q \geqslant c_1 n/\log n$ for some constants $c_0, c_1$, one has $a_1 = O(n/q) = O(\log n)$, $a_2(3) = o(n^2)$. Assumption (A2) yields that $\alpha([\frac{n}{q+1}])^{6/7} \leq Cn^{-6\lambda_0 c_0/7}$. Thus, for $n$ large enough,

$$(A.4) \quad P\left\{\frac{1}{n}\left|\sum_{l=1}^n \omega_J^*(\mathbf{X}_l, x_1)\right| > \frac{\rho \log n}{\sqrt{nh}}\right\} \leq cn^{-c_2\rho^2}\log n + Cn^{2-6\lambda_0 c_0/7}.$$

We divide the interval $[0, 1]$ into $M_n \sim n^6$ equally spaced intervals with disjoint endpoints $0 = x_{1,0} < x_{1,1} < \cdots < x_{1,M_n} = 1$. Employing the discretization method,

$$\sup_{x_1 \in [0,1]} \sup_{1 \leq J \leq N} \left|n^{-1}\sum_{l=1}^n \omega_J^*(\mathbf{X}_l, x_1)\right|$$

$$(A.5) \quad = \sup_{0 \leq k \leq M_n} \sup_{1 \leq J \leq N} \left|n^{-1}\sum_{l=1}^n \omega_J^*(\mathbf{X}_l, x_{1,k})\right|$$

$$+ \sup_{1 \leq k \leq M_n} \sup_{1 \leq J \leq N} \sup_{x_1 \in [x_{1,k-1}, x_{1,k}]} \left|n^{-1}\sum_{l=1}^n \{\omega_J^*(\mathbf{X}_l, x_1) - \omega_J^*(\mathbf{X}_l, x_{1,k})\}\right|.$$

By (A.4), there exists a large enough value $\rho > 0$ such that for any $1 \leq k \leq M_n$,

$$P\left\{\frac{1}{n}\left|\sum_{l=1}^{n}\omega_J^*(\mathbf{X}_l, x_{1,k})\right| > \rho(nh)^{-1/2}\log n\right\} \leq n^{-10}, \qquad 1 \leq J \leq N,$$

which implies that

$$\sum_{n=1}^{\infty} P\left\{\sup_{0 \leq k \leq M_n} \sup_{1 \leq J \leq N} \left|n^{-1}\sum_{l=1}^{n}\omega_J^*(\mathbf{X}_l, x_{1,k})\right| \geq \rho\frac{\log n}{\sqrt{nh}}\right\}$$

$$\leq \sum_{n=1}^{\infty}\sum_{k=1}^{M_n}\sum_{J=1}^{N} P\left\{\left|n^{-1}\sum_{l=1}^{n}\omega_J^*(\mathbf{X}_l, x_{1,k})\right| \geq \rho\frac{\log n}{\sqrt{nh}}\right\}$$

$$\leq \sum_{n=1}^{\infty} NM_n n^{-10} < \infty.$$

Thus, the Borel–Cantelli lemma entails that

(A.6) $$\sup_{0 \leq k \leq M_n} \sup_{1 \leq J \leq N} \left|n^{-1}\sum_{l=1}^{n}\omega_J^*(\mathbf{X}_l, x_{1,k})\right| = O_p(\log n/\sqrt{nh}).$$

Employing Lipschitz continuity of the kernel $K$, for $x_1 \in [x_{1,k-1}, x_{1,k}]$

$$\sup_{1 \leq k \leq M_n} |K_h(X_{l1} - x_1) - K_h(X_{l1} - x_{1,k})| \leq C_K M_n^{-1} h^{-2}.$$

According to the fact that $M_n \sim n^6$, one has

$$\sup_{1 \leq k \leq M_n} \sup_{1 \leq J \leq N} \sup_{x_1 \in [x_{1,k-1}, x_{1,k}]} \left|n^{-1}\sum_{l=1}^{n}\{\omega_J^*(\mathbf{X}_l, x_1) - \omega_J^*(\mathbf{X}_l, x_{1,k})\}\right|$$

$$\leq C_K M_n^{-1} h^{-2} \sup_{x_2 \in [0,1]} \sup_{1 \leq J \leq N} |B_{J,2}(x_2)|$$

$$= O(M_n^{-1} h^{-2} H^{-1/2}) = o(n^{-1}).$$

Thus (A.2) follows instantly from (A.5) and (A.6). As a result of Lemma A.3 and (A.2), (A.3) holds. □

LEMMA A.5. *Under assumptions* (A4) *and* (A6), *there exist constants* $C_0 > c_0 > 0$ *such that for any* $\mathbf{a} = (a_0, a_{1,1}, \ldots, a_{N,1}, a_{1,2}, \ldots, a_{N,2})$,

(A.7) $$c_0\left(a_0^2 + \sum_{J,\alpha} a_{J,\alpha}^2\right) \leq \left\|a_0 + \sum_{J,\alpha} a_{J,\alpha} B_{J,\alpha}\right\|_2^2 \leq C_0\left(a_0^2 + \sum_{J,\alpha} a_{J,\alpha}^2\right).$$

We refer the proof of the above lemma to Lemma A.4 in [26]. The next lemma provides the size of $\tilde{\mathbf{a}}^T\tilde{\mathbf{a}}$, where $\tilde{\mathbf{a}}$ is the least squares solution defined by (3.6).

LEMMA A.6. *Under assumptions* (A2)–(A6), *$\tilde{\mathbf{a}}$ satisfies*

$$(A.8) \qquad \tilde{\mathbf{a}}^T\tilde{\mathbf{a}} = \tilde{a}_0^2 + \sum_{J=1}^{N}\sum_{\alpha=1}^{2}\tilde{a}_{J,\alpha}^2 = O_p\{N(\log n)^2/n\}.$$

PROOF. According to (3.7) and (3.8), $\tilde{\mathbf{a}}^T\mathbf{B}^T\mathbf{B}\tilde{\mathbf{a}} = \tilde{\mathbf{a}}^T(\mathbf{B}^T\mathbf{E})$. Thus

$$(A.9) \qquad \|\mathbf{B}\tilde{\mathbf{a}}\|_{2,n}^2 = \tilde{\mathbf{a}}^T\begin{pmatrix} \mathbf{1} & \\ & \langle B_{J,\alpha}, B_{J',\alpha'}\rangle_{2,n} \end{pmatrix}\tilde{\mathbf{a}} = \tilde{\mathbf{a}}^T(n^{-1}\mathbf{B}^T\mathbf{E}).$$

By (A.15), $\|\mathbf{B}\tilde{\mathbf{a}}\|_{2,n}^2$ is bounded below in probability by $(1-A_n)\|\mathbf{B}\tilde{\mathbf{a}}\|_2^2$. According to (A.7), one has

$$(A.10) \qquad \|\mathbf{B}\tilde{\mathbf{a}}\|_2^2 = \left\|\tilde{a}_0^2 + \sum_{J=1}^{N}\sum_{\alpha=1}^{2}\tilde{a}_{J,\alpha}^2\right\|_2^2 \geq c_0\left(\tilde{a}_0^2 + \sum_{J,\alpha}\tilde{a}_{J,\alpha}^2\right).$$

Meanwhile one can show that $\tilde{\mathbf{a}}^T(n^{-1}\mathbf{B}^T\mathbf{E})$ is bounded above by

$$(A.11) \qquad \sqrt{\tilde{a}_0^2 + \sum_{J,\alpha}\tilde{a}_{J,\alpha}^2}\left[\left\{\frac{1}{n}\sum_{i=1}^{n}\sigma(\mathbf{X}_i)\varepsilon_i\right\}^2 + \sum_{J,\alpha}\left\{\frac{1}{n}\sum_{i=1}^{n}B_{J,\alpha}(X_{i\alpha})\sigma(\mathbf{X}_i)\varepsilon_i\right\}^2\right]^{1/2}.$$

Combining (A.9), (A.10) and (A.12), the squared norm $\tilde{\mathbf{a}}^T\tilde{\mathbf{a}}$ is bounded by

$$c_0^{-2}(1-A_n)^{-2}\left[\left\{\frac{1}{n}\sum_{i=1}^{n}\sigma(\mathbf{X}_i)\varepsilon_i\right\}^2 + \sum_{J,\alpha}\left\{\frac{1}{n}\sum_{i=1}^{n}B_{J,\alpha}(X_{i\alpha})\sigma(\mathbf{X}_i)\varepsilon_i\right\}^2\right].$$

Using the same truncation of $\varepsilon$ as in Lemma A.11, the Bernstein inequality entails that

$$\left|n^{-1}\sum_{i=1}^{n}\sigma(\mathbf{X}_i)\varepsilon_i\right| + \max_{1\leq J\leq N, \alpha=1,2}\left|n^{-1}\sum_{i=1}^{n}B_{J,\alpha}(X_{i\alpha})\sigma(\mathbf{X}_i)\varepsilon_i\right| = O_p(\log n/\sqrt{n}).$$

Therefore (A.8) holds since $A_n$ is of order $o_p(1)$. □

### A.2. Empirical approximation of the theoretical inner product.

LEMMA A.7. *Under assumptions (*A2)*, (*A4) *and* (A6)*, one has*

$$\sup_{J,\alpha}|\langle 1, B_{J,\alpha}\rangle_{2,n} - \langle 1, B_{J,\alpha}\rangle_2| = O_p(n^{-1/2}\log n), \tag{A.12}$$

$$\sup_{J,J',\alpha}|\langle B_{J,\alpha}, B_{J',\alpha}\rangle_{2,n} - \langle B_{J,\alpha}, B_{J',\alpha}\rangle_2| = O_p(n^{-1/2}H^{-1/2}\log n), \tag{A.13}$$

$$\sup_{1\leq J,J'\leq N, \alpha\neq\alpha'}|\langle B_{J,\alpha}, B_{J',\alpha'}\rangle_{2,n} - \langle B_{J,\alpha}, B_{J',\alpha'}\rangle_2|$$
$$= O_p(n^{-1/2}\log n). \tag{A.14}$$

We refer the proof of the above lemma to Lemma A.7 in [26].

LEMMA A.8. *Under assumptions* (A2), (A4) *and* (A6), *one has*

$$A_n = \sup_{g_1,g_2\in G^{(-1)}} \frac{|\langle g_1, g_2\rangle_{2,n} - \langle g_1, g_2\rangle_2|}{\|g_1\|_2\|g_2\|_2} = O_p\left(\frac{\log n}{n^{1/2}H^{1/2}}\right) = o_p(1). \tag{A.15}$$

PROOF. For every $g_1, g_2 \in G^{(-1)}$, one can write

$$g_1(X_1, X_2) = a_0 + \sum_{J=1}^{N}\sum_{\alpha=1}^{2} a_{J,\alpha}B_{J,\alpha}(X_\alpha),$$

$$g_2(X_1, X_2) = a_0' + \sum_{J'=1}^{N}\sum_{\alpha'=1}^{2} a_{J',\alpha'}'B_{J',\alpha'}(X_{\alpha'}),$$

where for any $J, J' = 1, \ldots, N, \alpha, \alpha' = 1, 2$, $a_{J,\alpha}$ and $a_{J',\alpha'}$ are real constants. Then

$$|\langle g_1, g_2\rangle_{2,n} - \langle g_1, g_2\rangle_2| \leq \left|\sum_{J,\alpha}\langle a_0', a_{J,\alpha}B_{J,\alpha}\rangle_{2,n}\right| + \left|\sum_{J',\alpha'}\langle a_0, a_{J',\alpha'}'B_{J',\alpha'}\rangle_{2,n}\right|$$
$$+ \sum_{J,J',\alpha,\alpha'}|a_{J,\alpha}||a_{J',\alpha'}'||\langle B_{J,\alpha}, B_{J',\alpha'}\rangle_{2,n}$$
$$- \langle B_{J,\alpha}, B_{J',\alpha'}\rangle_2|$$
$$= L_1 + L_2 + L_3.$$

The equivalence of norms given in (A.7) and (A.12) leads to

$$L_1 \leq A_{n,1}\cdot|a_0'|\cdot\sum_{J,\alpha}|a_{J,\alpha}|$$

$$\leq C_0 A_{n,1} \left( a_0'^2 + \sum_{J,\alpha} a_{J,\alpha}'^2 \right)^{1/2} \left( \sum_{J,\alpha} a_{J,\alpha}^2 \right)^{1/2} N^{1/2}$$

$$\leq C_{A,1} A_{n,1} \|g_1\|_2 \|g_2\|_2 H^{-1/2}$$

$$= O_p(n^{-1/2} H^{-1/2} \log n) \|g_1\|_2 \|g_2\|_2.$$

Similarly, $L_2 = O_p(n^{-1/2} H^{-1/2} \log n) \|g_1\|_2 \|g_2\|_2$. By the Cauchy–Schwarz inequality

$$L_3 \leq \sum_{J,J',\alpha,\alpha'} |a_{J,\alpha}| |a'_{J',\alpha'}| \max(A_{n,2}, A_{n,3})$$

$$\leq C_{A,2} \max(A_{n,2}, A_{n,3}) \|g_1\|_2 \|g_2\|_2$$

$$= O_p(n^{-1/2} H^{-1/2} \log n) \|g_1\|_2 \|g_2\|_2.$$

Therefore, statement (A.15) is established. □

**A.3. Proof of Lemma 5.2.** Denote $\mathbf{V}$ as the theoretical inner product of the B spline basis $\{1, B_{J,\alpha}(x_\alpha), J = 1, \ldots, N, \alpha = 1, 2\}$, that is,

(A.16) $$\mathbf{V} = \begin{pmatrix} 1 & \mathbf{0}_{2N}^T \\ \mathbf{0}_{2N} & \langle B_{J,\alpha}, B_{J',\alpha'} \rangle_2 \end{pmatrix}_{\substack{1 \leq \alpha, \alpha' \leq 2, \\ 1 \leq J, J' \leq N}},$$

where $\mathbf{0}_p = \{0, \ldots, 0\}^T$. Let $\mathbf{S}$ be the inverse matrix of $\mathbf{V}$, that is,

(A.17) $$\mathbf{S} = \mathbf{V}^{-1} = \begin{pmatrix} 1 & \mathbf{0}_N^T & \mathbf{0}_N^T \\ \mathbf{0}_N & \mathbf{V}_{11} & \mathbf{V}_{12} \\ \mathbf{0}_N & \mathbf{V}_{21} & \mathbf{V}_{22} \end{pmatrix}^{-1} = \begin{pmatrix} 1 & \mathbf{0}_N^T & \mathbf{0}_N^T \\ \mathbf{0}_N & \mathbf{S}_{11} & \mathbf{S}_{12} \\ \mathbf{0}_N & \mathbf{S}_{21} & \mathbf{S}_{22} \end{pmatrix}.$$

LEMMA A.9. *Under assumptions* (A4) *and* (A6), *for* $\mathbf{V}$, $\mathbf{S}$ *defined in* (A.16), (A.17), *there exist constants* $C_V > c_V > 0$ *and* $C_S > c_S > 0$ *such that*

$$c_V \mathbf{I}_{2N+1} \leq \mathbf{V} \leq C_V \mathbf{I}_{2N+1}, \qquad c_S \mathbf{I}_{2N+1} \leq \mathbf{S} \leq C_S \mathbf{I}_{2N+1}.$$

We refer the proof of the above lemma to Lemma A.9 in [26]. Next we denote

$$\mathbf{V}^* = \begin{pmatrix} 0 & \mathbf{0}_{2N}^T \\ \mathbf{0}_{2N} & \langle B_{J,\alpha}, B_{J',\alpha'} \rangle_{2,n} - \langle B_{J,\alpha}, B_{J',\alpha'} \rangle_2 \end{pmatrix}_{\substack{1 \leq \alpha, \alpha' \leq 2. \\ 1 \leq J, J' \leq N}}$$

Then $\tilde{\mathbf{a}}$ in (3.8) can be rewritten as

(A.18) $$\tilde{\mathbf{a}} = (\mathbf{B}^T \mathbf{B})^{-1} \mathbf{B}^T \mathbf{E} = \left( \frac{1}{n} \mathbf{B}^T \mathbf{B} \right)^{-1} \left( \frac{1}{n} \mathbf{B}^T \mathbf{E} \right)$$

$$= (\mathbf{V} + \mathbf{V}^*)^{-1} \left( \frac{1}{n} \mathbf{B}^T \mathbf{E} \right).$$

Now define $\hat{\mathbf{a}} = \{\hat{a}_0, \hat{a}_{1,1}, \ldots, \hat{a}_{N,1}, \hat{a}_{1,2}, \ldots, \hat{a}_{N,2}\}^T$ as

(A.19) $$\hat{\mathbf{a}} = \mathbf{V}^{-1}(n^{-1}\mathbf{B}^T\mathbf{E}) = \mathbf{S}(n^{-1}\mathbf{B}^T\mathbf{E}),$$

and define a theoretical version of $\Psi_v^{(2)}(x_1)$ in (5.6) as

(A.20) $$\hat{\Psi}_v^{(2)}(x_1) = n^{-1}\sum_{i=1}^{n}\sum_{J=1}^{N}\hat{a}_{J,2}\omega_J(\mathbf{X}_i, x_1).$$

LEMMA A.10. *Under assumptions* (A2) *to* (A6),
$$\sup_{x_1 \in [0,1]} |\Psi_v^{(2)}(x_1) - \hat{\Psi}_v^{(2)}(x_1)| = O_p\{(\log n)^2/(nH)\}.$$

PROOF. According to (A.18) and (A.19), one has $\mathbf{V}\hat{\mathbf{a}} = (\mathbf{V} + \mathbf{V}^*)\tilde{\mathbf{a}}$, which implies that $\mathbf{V}^*\tilde{\mathbf{a}} = \mathbf{V}(\hat{\mathbf{a}} - \tilde{\mathbf{a}})$. Using (A.13) and (A.14), one obtains that

$$\|\mathbf{V}(\hat{\mathbf{a}} - \tilde{\mathbf{a}})\| = \|\mathbf{V}^*\tilde{\mathbf{a}}\| \leq O_p(n^{-1/2}H^{-1}\log n)\|\tilde{\mathbf{a}}\|.$$

According to Lemma A.6, $\|\tilde{\mathbf{a}}\| = O_p(n^{-1/2}N^{1/2}\log n)$, so one has

$$\|\mathbf{V}(\hat{\mathbf{a}} - \tilde{\mathbf{a}})\| \leq O_p\{(\log n)^2 n^{-1}N^{3/2}\}.$$

By Lemma A.9, $\|(\hat{\mathbf{a}} - \tilde{\mathbf{a}})\| = O_p\{(\log n)^2 n^{-1}N^{3/2}\}$. Lemma A.6 then implies

(A.21) $$\|\hat{\mathbf{a}}\| \leq \|(\hat{\mathbf{a}} - \tilde{\mathbf{a}})\| + \|\tilde{\mathbf{a}}\| = O_p(\log n\sqrt{N/n}).$$

Additionally, $|\Psi_v^{(2)}(x_1) - \hat{\Psi}_v^{(2)}(x_1)| = |\sum_{J=1}^{N}(\tilde{a}_{J,2} - \hat{a}_{J,2})\frac{1}{n}\sum_{l=1}^{n}\omega_J(\mathbf{X}_l, x_1)|$. So

$$\sup_{x \in [0,1]} |\Psi_v^{(2)}(x_1) - \hat{\Psi}_v^{(2)}(x_1)| \leq \sqrt{N}O_p\left\{\frac{(\log n)^2}{nH}\right\}O_p(H^{1/2}) = O_p\left\{\frac{(\log n)^2}{nH}\right\}.$$

Therefore the lemma follows. □

LEMMA A.11. *Under assumptions* (A2)–(A6), *for* $\hat{\Psi}_v^{(2)}(x_1)$ *as defined in* (A.20), *one has*

$$\sup_{x_1 \in [0,1]} |\hat{\Psi}_v^{(2)}(x_1)| = \sup_{x_1 \in [0,1]}\left|n^{-1}\sum_{i=1}^{n}K_h(X_{i1} - x_1)\sum_{J=1}^{N}\hat{a}_{J,2}B_{J,2}(X_{i2})\right| = O_p(H).$$

PROOF. Note that

(A.22) $$|\hat{\Psi}_v^{(2)}(x_1)| \leq \left|\sum_{J=1}^{N}\hat{a}_{J,2}\mu_{\omega_J}(x_1)\right|$$
$$+ \left|\sum_{J=1}^{N}\hat{a}_{J,2}n^{-1}\sum_{i=1}^{n}\{\omega_J(\mathbf{X}_i, x_1) - \mu_{\omega_J}(x_1)\}\right|$$
$$= Q_1(x_1) + Q_2(x_1).$$

By the Cauchy–Schwarz inequality, (A.21) Lemma A.4 and assumptions (A5), (A6),

$$(\text{A.23}) \quad \sup_{x_1 \in [0,1]} Q_2(x_1) = O_p(\log n \sqrt{N/n}) \sqrt{N} O_p\left(\frac{\log n}{\sqrt{nh}}\right) = O_p\left\{\frac{(\log n)^3}{\sqrt{n}}\right\}.$$

Using the discretization idea again as in the proof of Lemma A.4, one has

$$\sup_{x_1 \in [0,1]} Q_1(x_1)$$

$$(\text{A.24}) \quad \leq \max_{1 \leq k \leq M_n} \left| \sum_{J=1}^N \hat{a}_{J,2} \mu_{\omega_J}(x_{1,k}) \right|$$

$$+ \max_{1 \leq k \leq M_n} \sup_{x_1 \in [x_{1,k-1}, x_{1,k}]} \left| \sum_{J=1}^N \hat{a}_{J,2} \mu_{\omega_J}(x_1) - \sum_{J=1}^N \hat{a}_{J,2} \mu_{\omega_J}(x_{1,k}) \right|$$

$$= T_1 + T_2,$$

where $M_n \sim n$. Define next

$$W_1 = \max_{1 \leq k \leq M_n} \left| n^{-1} \sum_{1 \leq i \leq n} \sum_{1 \leq J, J' \leq N} \mu_{\omega_J}(x_{1,k}) s_{J+N+1, J'+1} B_{J',1}(X_{i1}) \sigma(\mathbf{X}_i) \varepsilon_i \right|,$$

$$W_2 = \max_{1 \leq k \leq M_n} \left| n^{-1} \sum_{1 \leq i \leq n} \sum_{1 \leq J, J' \leq N} \mu_{\omega_J}(x_{1,k}) s_{J+N+1, J'+N+1} B_{J',2}(X_{i2}) \sigma(\mathbf{X}_i) \varepsilon_i \right|.$$

Then it is clear that $T_1 \leq W_1 + W_2$. Next we will show that $W_1 = O_p(H)$. Let $D_n = n^{\theta_0}(\frac{1}{2+\delta} < \theta_0 < \frac{2}{5})$, where $\delta$ is the same as in assumption (A3). Define

$$\varepsilon_{i,D}^- = \varepsilon_i I(|\varepsilon_i| \leq D_n), \qquad \varepsilon_{i,D}^+ = \varepsilon_i I(|\varepsilon_i| > D_n), \qquad \varepsilon_{i,D}^* = \varepsilon_{i,D}^- - E(\varepsilon_{i,D}^- | \mathbf{X}_i),$$

$$U_{i,k} = \mu_\omega(x_{1,k})^T \mathbf{S}_{21} \{B_{1,1}(X_{i1}), \ldots, B_{1,N}(X_{i1})\}^T \sigma(\mathbf{X}_i) \varepsilon_{i,D}^*.$$

Denote $W_1^D = \max_{1 \leq k \leq M_n} |n^{-1} \sum_{i=1}^n U_{i,k}|$ as the truncated centered version of $W_1$. Next we show that $|W_1 - W_1^D| = O_p(H)$. Note that $|W_1 - W_1^D| \leq \Lambda_1 + \Lambda_2$, where

$$\Lambda_1 = \max_{1 \leq k \leq M_n} \left| \frac{1}{n} \sum_{i=1}^n \sum_{1 \leq J, J' \leq N} \mu_{\omega_J}(x_{1,k}) s_{J+N+1, J'+1} \right.$$

$$\left. \times B_{J',1}(X_{i1}) \sigma(\mathbf{X}_i) E(\varepsilon_{i,D}^- | \mathbf{X}_i) \right|,$$

$$\Lambda_2 = \max_{1 \leq k \leq M_n} \left| \frac{1}{n} \sum_{i=1}^n \sum_{1 \leq J, J' \leq N} \mu_{\omega_J}(x_{1,k}) s_{J+N+1, J'+1} B_{J',1}(X_{i1}) \sigma(\mathbf{X}_i) \varepsilon_{i,D}^+ \right|.$$

Let $\mu_\omega(x_{1,k}) = \{\mu_{\omega_1}(x_{1,k}), \ldots, \mu_{\omega_N}(x_{1,k})\}^T$; then

$$\Lambda_1 = \max_{1 \leq k \leq M_n} \left| \mu_\omega(x_{1,k})^T \mathbf{S}_{21} \left\{ n^{-1} \sum_{i=1}^n B_{J',1}(X_{i1}) \sigma(\mathbf{X}_i) E(\varepsilon_{i,D}^- | \mathbf{X}_i) \right\}_{J'=1}^N \right|$$

$$\leq C_S \max_{1 \leq k \leq M_n} \left\{ \sum_{J=1}^N \mu_{\omega_J}^2(x_{1,k}) \sum_{J=1}^N \left\{ \frac{1}{n} \sum_{i=1}^n B_{J,1}(X_{i1}) \sigma(\tilde{\mathbf{X}}_i) E(\varepsilon_{i,D}^- | \mathbf{X}_i) \right\}^2 \right\}^{1/2}.$$

By assumption (A3), one has $|E(\varepsilon_{i,D}^- | \mathbf{X}_i)| = |E(\varepsilon_{i,D}^+ | \mathbf{X}_i)| \leq M_\delta D_n^{-(1+\delta)}$ and Lemma A.1 entails that $\sup_{J,\alpha} |\frac{1}{n} \sum_{i=1}^n B_{J,1}(X_{i1}) \sigma(\mathbf{X}_i)| = O_p(\log n/\sqrt{n})$. Therefore

$$\Lambda_1 \leq M_\delta D_n^{-(1+\delta)} \max_{1 \leq k \leq M_n} \left[ \sum_{J=1}^N \mu_{\omega_J}^2(x_{1,k}) \sum_{J=1}^N \left\{ \frac{1}{n} \sum_{i=1}^n B_{J,1}(X_{i1}) \sigma(\mathbf{X}_i) \right\}^2 \right]^{1/2}$$

$$= O_p\{N D_n^{-(1+\delta)} \log^2 n / n\} = O_p(H),$$

where the last step follows from the choice of $D_n$. Meanwhile

$$\sum_{n=1}^\infty P(|\varepsilon_n| \geq D_n) \leq \sum_{n=1}^\infty \frac{E|\varepsilon_n|^{2+\delta}}{D_n^{2+\delta}} = \sum_{n=1}^\infty \frac{E(E|\varepsilon_n|^{2+\delta} | \mathbf{X}_n)}{D_n^{2+\delta}} \leq \sum_{n=1}^\infty \frac{M_\delta}{D_n^{2+\delta}} < \infty,$$

since $\delta > 1/2$. By the Borel–Cantelli lemma, one has with probability 1,

$$n^{-1} \sum_{i=1}^n \sum_{1 \leq J, J' \leq N} \mu_{\omega_J}(x_{1,k}) s_{J+N+1, J'+1} B_{J',1}(X_{i1}) \sigma(\mathbf{X}_i) \varepsilon_{i,D}^+ = 0,$$

for large $n$. Therefore, one has $|W_1 - W_1^D| \leq \Lambda_1 + \Lambda_2 = O_p(H)$. Next we will show that $W_1^D = O_p(H)$. Note that the variance of $U_{i,k}$ is

$$\mu_\omega(x_{1,k})^T \mathbf{S}_{21} \text{var}(\{B_{1,1}(X_{i1}), \ldots, B_{N,1}(X_{i1})\}^T \sigma(\mathbf{X}_i) \varepsilon_{i,D}^*) \mathbf{S}_{21} \mu_\omega(x_{1,k}).$$

By assumption (A3), $c_\sigma^2 \mathbf{V}_{11} \leq \text{var}(\{B_{1,1}(X_{i1}), \ldots, B_{N,1}(X_{i1})\}^T \sigma(\mathbf{X}_i)) \leq C_\sigma^2 \mathbf{V}_{11}$,

$$\text{var}(U_{i,k}) \sim \mu_\omega(x_{1,k})^T \mathbf{S}_{21} \mathbf{V}_{11} \mathbf{S}_{21} \mu_\omega(x_{1,k}) V_{\varepsilon,D} = \mu_\omega(x_{1,k})^T \mathbf{S}_{21} \mu_\omega(x_{1,k}) V_{\varepsilon,D},$$

where $V_{\varepsilon,D} = \text{var}\{\varepsilon_{i,D}^* | \mathbf{X}_i\}$. Let $\kappa(x_{1,k}) = \{\mu_\omega(x_{1,k})^T \mu_\omega(x_{1,k})\}^{1/2}$. Then

$$c_S c_\sigma^2 \{\kappa(x_{1,k})\}^2 V_{\varepsilon,D} \leq \text{var}(U_{i,k}) \leq C_S C_\sigma^2 \{\kappa(x_{1,k})\}^2 V_{\varepsilon,D}.$$

Simple calculation leads to

$$E|U_{i,k}|^r \leq \{c_0 \kappa(x_{1,k}) D_n H^{-1/2}\}^{r-2} r! E|U_{i,k}|^2 < +\infty$$

for $r \geq 3$, so $\{U_{i,k}\}_{i=1}^n$ satisfies Cramér's condition with Cramér's constant $c_* = c_0 \kappa(x_{1,k}) D_n H^{-1/2}$; hence by Bernstein's inequality

$$P\left\{ \left| n^{-1} \sum_{l=1}^n U_{i,k} \right| \geq \rho_n \right\} \leq a_1 \exp\left( -\frac{q \rho_n^2}{25 m_2^2 + 5 c_* \rho_n} \right) + a_2(3) \alpha\left( \left[ \frac{n}{q+1} \right] \right)^{6/7},$$

where $m_2^2 \sim \{\kappa(x_{1,k})\}^2 V_{\varepsilon,D}, m_3 \leq \{c\{\kappa(x_{1,k})\}^3 H^{-1/2} D_n V_{\varepsilon,D}\}^{1/3}$,

$$\rho_n = \rho H, \qquad a_1 = 2\frac{n}{q} + 2\left(1 + \frac{\rho_n^2}{25m_2^2 + 5c_*\rho_n}\right), \qquad a_2(3) = 11n\left(1 + \frac{5m_3^{6/7}}{\rho_n}\right).$$

Similar arguments as in Lemma A.4 yield that as $n \to \infty$

$$\frac{q\rho_n^2}{25m_2^2 + 5c_*\rho_n} \sim \frac{q\rho_n}{c_*} = \frac{\rho n^{2/5}}{c_0(\log n)^{5/2} D_n} \to +\infty.$$

Taking $c_0, \rho$ large enough, one has

$$P\left\{\frac{1}{n}\left|\sum_{i=1}^n U_{i,k}\right| > \rho H\right\} \leq c\log n \exp\{-c_2\rho^2 \log n\} + Cn^{2-6\lambda_0 c_0/7} \leq n^{-3},$$

for $n$ large enough. Hence

$$\sum_{n=1}^\infty P(|W_1^D| \geq \rho H) = \sum_{n=1}^\infty \sum_{k=1}^{M_n} P\left(\left|\frac{1}{n}\sum_{i=1}^n U_{i,k}\right| \geq \rho H\right) \leq \sum_{n=1}^\infty M_n n^{-3} < \infty.$$

Thus, the Borel–Cantelli lemma entails that $W_1^D = O_p(H)$. Noting the fact that $|W_1 - W_1^D| = O_p(H)$, one has that $W_1 = O_p(H)$. Similarly $W_2 = O_p(H)$. Thus

(A.25) $$T_1 \leq W_1 + W_2 = O_p(H).$$

Employing the Cauchy–Schwarz inequality and Lipschitz continuity of the kernel $K$, assumption (A5), Lemma A.2(ii) and (A.21) lead to

(A.26) $$T_2 \leq O_p\left(\frac{N^{1/2}\log n}{n^{1/2}}\right)\frac{\{\sum_{J=1}^N EB_{J,2}^2(X_{12})\}^{1/2}}{h^2 M_n} = o_p(n^{-1/2}).$$

Combining (A.24), (A.25) and (A.26), one has $\sup_{x_1 \in [0,1]} Q_1(x_1) = O_p(H)$. The desired result follows from (A.23) and (A.23). $\square$

**Acknowledgments.** This work is part of the first author's dissertation under the supervision of the second author. The authors are very grateful to the Editor, Jianqing Fan, and three anonymous referees for their helpful comments.

ADDITIVE AUTOREGRESSION MODEL 31[26] WANG, L. and YANG, L. (2006). Spline-backfitted kernel smoothing of nonlinear additive autoregression model. Manuscript. Available at www.arxiv.org/abs/math/0612677.
[27] XUE, L. and YANG, L. (2006). Estimation of semiparametric additive coefficient model. *J. Statist. Plann. Inference* **136** 2506–2534. MR2279819
[28] XUE, L. and YANG, L. (2006). Additive coefficient modeling via polynomial spline. *Statist. Sinica* **16** 1423–1446. MR2327498
[29] YANG, L., HÄRDLE, W. and NIELSEN, J. P. (1999). Nonparametric autoregression with multiplicative volatility and additive mean. *J. Time Ser. Anal.* **20** 579–604. MR1720162
[30] YANG, L., SPERLICH, S. and HÄRDLE, W. (2003). Derivative estimation and testing in generalized additive models. *J. Statist. Plann. Inference* **115** 521–542. MR1985882
DEPARTMENT OF STATISTICS
UNIVERSITY OF GEORGIA
ATHENS, GEORGIA 30602
USA
E-MAIL: lilywang@uga.edu

DEPARTMENT OF STATISTICS
AND PROBABILITY
MICHIGAN STATE UNIVERSITY
EAST LANSING, MICHIGAN 48824
USA
E-MAIL: yangli@stt.msu.edu